\documentclass{amsproc}

\usepackage{amsmath,amssymb}
\newtheorem{theorem}{Theorem}[section]
\newtheorem{lemm}[theorem]{Lemma}
\newtheorem{prop}[theorem]{Proposition}

\theoremstyle{definition}
\newtheorem{defi}[theorem]{Definition}

\newtheorem{coro}[theorem]{Corollary}
\theoremstyle{remark}
\newtheorem{remark}[theorem]{Remark}

\numberwithin{equation}{section}

\def\frak{\mathfrak}

\def\cal{\mathcal}

\def\om{\omega}
\def\lra{\longrightarrow}

\def\dim{\hbox{dim}}

\def\om{\omega}

\def\mod{\hbox{mod}}

\newfont{\df}{eufm10}
\def\vep{\varepsilon}

\def\ot{\otimes}

\def\dim{\hbox{\rm dim}\,}

\def\ot{\otimes}

\def\ra{\rangle}
\def\la{\langle}

\begin{document}

\title[Two-parameter Quantum Affine Algebra, Drinfel'd Realization]
{Two-parameter Quantum Affine Algebra $U_{r,s}(\widehat{\frak
{sl}_n})$, Drinfel'd Realization and Quantum Affine Lyndon Basis}

\author[Hu]{Naihong Hu}
\address{Department of Mathematics, East China Normal University,
Shanghai 200062, PR China}
\email{nhhu@math.ecnu.edu.cn}

\author[Rosso]{Marc Rosso}
\address{D\'epartment math\'ematiques et applications, Ecole Normale Superieure, 45 Rue de Ulm, 75230 Paris Cedex 05, France}
\email{Marc.Rosso@ens.fr}

\author[Zhang]{Honglian Zhang}
\address{Department of Mathematics,
East China Normal University, Shanghai 200062, PR
China}\curraddr{Department of Mathematics, Shanghai University,
Shanghai 200444, PR China} \email{hlzhangmath@shu.edu.cn}
\thanks{N.H.,
supported in part by the NNSF (Grant 10431040), the TRAPOYT and the
FUDP from the MOE of China, the SRSTP from the STCSM, die Deutche
Forschungsgemeinschaft (DFG), as well as an ICTP long-term visiting
scholarship. H.Z., supported by a Ph.D. Program Scholarship Fund of
ECNU 2006. }

\subjclass{Primary 17B37, 81R50; Secondary 17B35}

\keywords{Quantum affine algebra, Drinfel'd realization, Drinfel'd
isomorphism, quantum ``affine" Lyndon basis.}
\begin{abstract}
We further define two-parameter quantum affine algebra
$U_{r,s}(\widehat{\frak {sl}_n})$ $(n>2)$ after the work on the
finite cases (see [BW1], [BGH1], [HS] \& [BH]), which turns out to
be a Drinfel'd double. Of importance for the quantum {\it affine}
cases is that we can work out the compatible two-parameter version
of the Drinfel'd realization as a quantum affinization of
$U_{r,s}(\frak{sl}_n)$ and establish the Drinfel'd isomorphism
Theorem in the two-parameter setting, via developing a new
combinatorial approach (quantum calculation) to the quantum {\it
affine} Lyndon basis we present  (with an explicit valid algorithm
based on the use of Drinfel'd generators).
\end{abstract}

\maketitle

\section{ Introduction}
\medskip

\noindent {\bf 1.1} \  In 2001, Benkart-Witherspoon investigated the
structures of two-parameter quantum groups $U_{r,s}(\frak g)$ for
${\frak g}={\frak{gl_n}}$, or ${\frak{sl_n}}$ in [BW1] originally
obtained by Takeuchi [T], and the finite-dimensional weight
representation theory in [BW2], and further obtained some new
finite-dimensional pointed Hopf algebras in [BW3] when $rs^{-1}$ is
a root of unity, which possess new ribbon elements under some
conditions (and will yield new invariants of knots and links). These
show that two-parameter quantum groups are well worth further
studying.

\medskip \noindent {\bf 1.2} \ In 2004, Bergeron-Gao-Hu [BGH1] gave
the structures of two-parameter quantum groups $U_{r,s}({\frak{g}})$
for ${\frak g}=\frak {so}_{2n+1},\, \frak{sp}_{2n},\,
\frak{so}_{2n}$, and developed in [BGH2] the highest weight
representation theory when $rs^{-1}$ is not a root of unity.
Especially, [BGH1] explored the environment condition upon which the
Lusztig's symmetries exist for the classical simple Lie algebras
$\frak g$, namely, they exist as $\mathbb Q$-isomorphisms between
$U_{r,s}(\frak g)$ and the associated object
$U_{s^{-1},r^{-1}}(\frak g)$ only when $\text{\rm rank}\,(\frak
g)=2$, and in the case when $\text{\rm rank}\,(\frak g)>2$, the
sufficient and necessary condition for the existence of Lusztig's
symmetries between $U_{r,s}(\frak g)$ and its associated object
forces $U_{r,s}(\frak g)$ to take the ``one-parameter" form
$U_{q,q^{-1}}(\frak g)$ where $r=s^{-1}=q$. In other words, when
$\text{\rm rank}\,(\frak g)>2$, the Lusztig's symmetries exist only
for the one-parameter quantum groups $U_{q, q^{-1}}(\frak g)$ as
$\mathbb Q(q)$-automorphisms (rather merely as $\mathbb
Q$-isomorphisms). In this case, these symmetries give rise to, with
respect to modulo some identification of group-like elements, the
usual Lusztig symmetries on quantum groups $U_q(\frak g)$ of
Drinfel'd-Jimbo type. The Lusztig symmetry property indicates that
there do exist the remarkable differences between the two-parameter
quantum groups in question and the one-parameter quantum groups of
Drinfel'd-Jimbo type. Afterwards, Hu-Shi [HS] and Bai-Hu [BH]
studied the two-parameter quantum groups for type $G_2$ and $E$
cases. Through the work, we found that the treatments in
two-parameter cases are frequently more subtle to follow
combinatorial approaches only, for instance, the description of the
convex PBW-type basis (cf. [BH]) has to appeal to the use of Lyndon
words (see [R2] and references therein) because there is no braid
group available in question.

\smallskip
Thereby so far, it seems desirable to extend these kind of the
two-parameter quantum groups in the Benkart-Witherspoon's sense in
finite cases to the affine cases. The present paper is aimed at this
purpose for the affine type $A^{(1)}_n$ ($n>1$) case. To this end,
we first give the defining structure of $U_{r,s}(\widehat{\frak
{sl}_n})$ ($n>2$) (whereas $U_{r,s}(\widehat{\frak {sl}_2})$ is
essentially isomorphic to $U_{q,q^{-1}}(\widehat{\frak {sl}_2})$ if
set $rs^{-1}=q^2$, which is not considered in the paper).

\medskip \noindent {\bf 1.3} \  As is
well-known, the importance of the Drinfel'd generators (in the
Drinfel'd realization) for quantum affine algebras is just like that
of the loop generators (in the loop realization) for affine
Kac-Moody algebras (see [Ga], [K]). Early in 1987, Drinfel'd [Dr2]
put forward his famous new (conjectural) realization of quantum
affine algebras $U_q(\widehat {\frak g})$ with $\frak g$ semisimple,
because he recognized that the study of finite dimensional
representations of $U_q(\widehat {\frak g})$ is made easier by the
use of this realization on the set of Drinfel'd generators, which is
called the Drinfel'd realization of $U_q(\widehat {\frak g})$ or the
Drinfel'd quantum affinization of $U_q(\frak g)$. Besides this, the
Drinfel'd realization also finds its main contribution to the
construction of vertex representations for quantum affine algebras
$U_q(\widehat {\frak g})$ (see [FJ], [J1], [DI2], etc.), as does the
loop realization in the vertex representation theory of affine
Kac-Moody algebras (see [K]). In 1993, Khoroshkin-Tolstoy [KT]
constructed the Drinfel'd realization for the untwisted types using
a Cartan-Weyl generators system with no proof. The first perfect
proof of the Drinfel'd isomorphsim only for the untwisted types was
given by Beck [B2] till 1994 making use of his extended braid group
actions, based on the work of Damiani [Da],
Levendorskii-Soibel'man-Stukopin [LSS] for the case
$U_q(\widehat{\frak {sl}_2})$. In 1998, Jing [J2] basically adopted
the inverse map suggested by Beck for the untwisted types (see the
final remark in [B2, Section 4]) and gave a combinatorial proof for
the Drinfel'd isomorphism for the untwisted types.

\smallskip
\noindent {\bf 1.4} \ In order to further explore and enrich the
structure and representation theory of the two-parameter quantum
affine algebras later on, the another main result of this paper is
to give the Drinfel'd realization of $U_{r,s}(\widehat{\frak
{sl}_n})$ ($n>2$). Its definition depends on the self-compatible
defining system (Definition 3.1), which in the two-parameter
setting, varies dramatically in comparison with the one-parameter
cases (see [Dr2], or [B2, Theorem 4.7]) and is nontrivial to match
up here and there the whole relations together. Indeed, to invent
the two-parameter version of Drinfel'd realization needs some
insights, e.g., from the antisymmetric point of view via the
$\mathbb{Q}$-algebra antiautomorphism $\tau$, based on some
information from the combinatorial description of the convex
PBW-type basis via the Lyndon words (see [R2], [BH], etc.), and
also, the proof of the Drinfel'd isomorphism in our case depends
completely on the combinatorial approach with specific techniques to
design those defining relations in order to fit the compatibilities
in the whole system. If the readers go with us into the details,
they will find how our quantum calculations (in somehow a bit
tedious) work well and necessarily for exactly verifying the
compatibilities of the defining system. The reason is that the
method we expanded, to some extent, essentially follows an approach
to a kind of description of the quantum ``affine" Lyndon basis.
Actually, we can construct explicitly all quantum real and imaginary
root vectors using this method (see Lemmas 4.7 \& 4.8, together with
Definition 3.9).

\medskip \noindent {\bf 1.5} \ The paper is organized as follows.
We first give the structure of two-parameter quantum affine algebra
$U_{r,s}(\widehat{\frak {sl}_n})$ ($n>2$) as Hopf algebra in Section
2. We prove that two-parameter quantum affine algebra
$U_{r,s}(\widehat{\frak {sl}_n})$ is characterized as Drinfel'd
double ${\cal D}({\widehat{\cal B}}, {\widehat{\cal B'}})$ of Hopf
subalgebras ${\widehat{\cal B}}$, ${\widehat{\cal B'}}$ with respect
to a skew-dual pairing. In Section 3, we explicitly describe the
two-parameter Drinfel'd quantum affinization of
$U_{r,s}(\frak{sl}_n)$ ($n>2$), that is, the Drinfel'd realization
in two-parameter case which is antisymmetric with respect to the
$\mathbb{Q}$-algebra antiautomorphism $\tau$. In the case when
$rs=1$, i.e., $r=s^{-1}=q$, our result modulo some identification
yields the usual Drinfel'd realization of quantum affine algebra
$U_q(\widehat{\frak {sl}_n})$ of Drinfel'd-Jimbo type (see [Dr2],
[B2], [DI1], [J2], etc.). Since Beck's extended braid group actions
approach is invalid for our case, we combine the Lyndon words
description ([R2]) with the quantum Lie bracket operation ([J2]) to
develop a combinatorial trick in the quantum affine case (we call it
quantum calculations), which can be utilized in the construction of
all the quantum root vectors (including real and imaginary ones), so
that we can formulate and prove the quantum ``{\it affine}" Lyndon
basis for the first time (in a more explicit form than that of [B1])
for ${\mathcal U}_{r,s}(\widetilde{\frak {n}}^\pm)$ based on the
Drinfel'd realization in Section 3, and further prove the Drinfel'd
isomorphism using our combinatorial algorithm in Section 4. In fact,
our proof also provides a concrete process of how to construct the
Drinfel'd generators using the Chevalley-Kac-Lusztig generators.

\section{Quantum Affine Algebra $U_{r,s}(\widehat{\frak {sl}_n})$ and Drinfel'd Double}
\medskip

\noindent {\bf 2.1} Let $\mathbb{K}=\mathbb{Q}(r,s)$ denote a
field of rational functions with two-parameters $r$, $s$ ($r\ne
\pm s$). Assume $\Phi$ is a finite root system of type $A_{n-1}$
with $\Pi$ a base of simple roots. Regard $\Phi$ as a subset of a
Euclidean space $E=\mathbb{R}^n$ with an inner product $( \, , )$.
Set $I=\{1,\cdots,n-1\}$, $I_0=\{0\}\cup I$. Let $\varepsilon_1,
\varepsilon_2, \cdots, \varepsilon_n$ denote an orthonormal basis
of $E$, then we can take
$\Pi=\{\alpha_i=\varepsilon_i-\varepsilon_{i+1}\mid i\in I \}$ and
$\Phi=\{\varepsilon_i-\varepsilon_j\mid i \neq j\in I\}$. Let
$\delta$ denote the primitive imaginary root of $\widehat{\frak
{sl}_n}$. Take $\alpha_0=\delta-(\varepsilon_1-\varepsilon_n)$,
then $\Pi'=\{\alpha_i\mid i\in I_0\}$ is a base of simple roots of
affine Lie algebra ${\widehat{\frak{sl}_n}}$.

Let $A=(a_{ij})$ $(i,\,j\in I_0)$ be a generalized Cartan matrix
associated to affine Lie algebra $\widehat{\frak {sl}_n}$. Let
${\frak{h}}$ be a vector space over $\mathbb{K}$ with a basis $\{
\,h_0,\, h_1,\, \cdots,\,h_{n-1},\,d\, \}$ and define the linear
action of $\alpha_i\, (i\in I_0)$ on ${\frak{h}}$ by
$$\alpha_i(h_j)=a_{ji},  \qquad \alpha_i(d)=\delta_{i,0}, \quad\textit{for }\ j\in I_0.$$

Let $Q=\mathbb{Z}\alpha_0+\cdots+\mathbb{Z}\alpha_{n-1}$ denote the
root lattice of $\widehat{\frak {sl}_n}$. The standard nondegenerate
symmetric bilinear form $(\cdot \,,\cdot)$ on ${\frak{h^*}}$
satisfies
$$(\alpha_i, \alpha_j)=a_{ij}, \qquad (\delta, \alpha_i)=(\delta, \delta)=0,\quad \forall\; i,\,j \in I_0.$$
\begin{defi}
Let $U=U_{r,s}(\widehat{\frak {sl}_n})$ $(n>2)$ be the unital
associative algebra over $\mathbb{K}$ generated by the elements
$e_j,\, f_j,\, \omega_j^{\pm 1},\, \omega_j'^{\,\pm 1}\, (j\in
I_0),\, \gamma^{\pm\frac{1}2},\,\gamma'^{\pm\frac{1}2},\,
D^{\pm1}, D'^{\,\pm1}$ (called the Chevalley-Kac-Lusztig
generators), satisfying the following relations:

\medskip
\noindent $(\textrm{A1})$ \
$\gamma^{\pm\frac{1}2},\,\gamma'^{\pm\frac{1}2}$ are central with
$\gamma=\om_\delta$, $\gamma'=\om'_\delta$, $\gamma\gamma'=rs$, such
that $\omega_i\,\omega_i^{-1}=\omega_i'\,\omega_i'^{\,-1}=1
=DD^{-1}=D'D'^{-1}$, and
\begin{equation*}
\begin{split}[\,\omega_i^{\pm 1},\omega_j^{\,\pm 1}\,]&=[\,\om_i^{\pm1},
D^{\pm1}\,]=[\,\om_j'^{\,\pm1}, D^{\pm1}\,] =[\,\om_i^{\pm1},
D'^{\pm1}\,]=0\\
&=[\,\omega_i^{\pm 1},\omega_j'^{\,\pm 1}\,]=[\,\om_j'^{\,\pm1},
D'^{\pm1}\,]=[D'^{\,\pm1}, D^{\pm1}]=[\,\omega_i'^{\pm
1},\omega_j'^{\,\pm 1}\,].
\end{split}
\end{equation*}
 $(\textrm{A2})$ \ For $\,i\in I_0$
and $j\in I$,
\begin{equation*}
\begin{array}{ll}
& D\,e_i\,D^{-1}=r^{\delta_{0i}}\,e_i,\qquad\qquad\qquad\qquad\;
D\,f_i\,D^{-1}=r^{-\delta_{0i}}\,f_i,\\
&\omega_j\,e_i\,\omega_j^{\,-1}=r^{(\varepsilon_j,
\alpha_i)}s^{(\varepsilon_{j+1}, \alpha_i )}\,e_i,\qquad\quad
\omega_j\,f_i\,\omega_j^{\,-1}=r^{-(\varepsilon_j,
\alpha_i)}s^{-(\varepsilon_{j+1},\alpha_i)}\,f_i,\\
&\omega_0\,e_i\,\omega_0^{\,-1}=r^{-(\varepsilon_{i+1},
\alpha_0)}s^{(\varepsilon_1, \alpha_i )}\,e_i,\qquad\
\omega_0\,f_i\,\omega_0^{\,-1}=r^{(\varepsilon_{i+1},
\alpha_0)}s^{-(\varepsilon_1, \alpha_i )}\,f_i .
\end{array}
\end{equation*}\\
$(\textrm{A3})$ \ For $\,i\in I_0$ and $j\in I$,
\begin{equation*}
\begin{array}{ll}
& D'\,e_i\,D'^{-1}=s^{\delta_{0i}}\,e_i,\qquad\qquad\qquad\quad\ \
D'\,f_i\,D'^{-1}=s^{-\delta_{0i}}\,f_i,\\
&\omega'_j\,e_i\,\omega'^{\,-1}_j=s^{(\varepsilon_j,
\alpha_i)}r^{(\varepsilon_{j+1}, \alpha_i)}\,e_i, \qquad\ \
\omega'_j\,f_i\,\omega'^{\,-1}_j=s^{-(\varepsilon_j,
\alpha_i)}r^{-(\varepsilon_{j+1},\alpha_i)}\,f_i,\\
&\omega'_0\,e_i\,\omega'^{\,-1}_0=s^{-(\varepsilon_{i+1},
\alpha_0)}r^{(\epsilon_1,\alpha_i)}\,e_i, \qquad
\omega'_0\,f_i\,\omega'^{\,-1}_0=s^{(\varepsilon_{i+1},
\alpha_0)}r^{-(\epsilon_1,\alpha_i)}\,f_i.
\end{array}
\end{equation*}\\
$(\textrm{A4})$ \ For $\,i,\, j\in I_0$, we have
 $$[\,e_i, f_j\,]=\frac{\delta_{ij}}{r-s}(\omega_i-\omega'_i).$$
$(\textrm{A5})$ \ For $\,i,\,j\in I_0$, but $(i,j)\notin
\{\,(0,n-1), \ (n-1,0)\,\}$ with $a_{ij}=0$, we have
 $$[\,e_i, e_j\,]=0=[\,f_i, f_j\,].$$
$(\textrm{A6})$ \ For $\,i\in I_0$, we have the $(r,s)$-Serre
relations:
\begin{gather*}
e_i^2e_{i+1}-(r+s)\,e_ie_{i+1}e_i+(rs)\,e_{i+1}e_i^2=0,\\
e_ie_{i+1}^2-(r+s)\,e_{i+1}e_ie_{i+1}+(rs)\,e_{i+1}^2e_i=0,\\
e_{n-1}^2e_0-(r+s)\,e_{n-1}e_0e_{n-1}+(rs)\,e_0e_{n-1}^2=0,\\
e_{n-1}e_0^2-(r+s)\,e_0e_{n-1}e_0+(rs)\,e_0^2e_{n-1}=0.
\end{gather*}
$(\textrm{A7})$ \ For $\,i\in I_0$, we have the
$(r,s)$-Serre relations:
\begin{gather*}
f_i^2f_{i+1}-(r^{-1}+s^{-1})\,f_if_{i+1}f_i+(r^{-1}s^{-1})\,f_{i+1}f_i^2=0,\\
f_if_{i+1}^2-(r^{-1}+s^{-1})\,f_{i+1}f_if_{i+1}+(r^{-1}s^{-1})\,f_{i+1}^2f_i=0,\\
f_{n-1}^2f_0-(r^{-1}+s^{-1})\,f_{n-1}f_0f_{n-1}+(r^{-1}s^{-1})\,f_0f_{n-1}^2=0,\\
f_{n-1}f_0^2-(r^{-1}+s^{-1})\,f_0f_{n-1}f_0+(r^{-1}s^{-1})\,f_0^2f_{n-1}=0.
\end{gather*}

$U_{r,s}(\widehat{\frak {sl}_n})$ is a Hopf algebra with the
coproduct $\Delta$, the counit $\vep$ and the antipode $S$ defined
below: for $i\in I_0$, we have
\begin{gather*}
\Delta(\gamma^{\pm\frac{1}2})=\gamma^{\pm\frac{1}2}\otimes
\gamma^{\pm\frac{1}2}, \qquad
\Delta(\gamma'^{\,\pm\frac{1}2})=\gamma'^{\,\pm\frac{1}2}\otimes
\gamma'^{\,\pm\frac{1}2}, \\
\Delta(D^{\pm1})=D^{\pm1}\otimes D^{\pm1},\qquad
\Delta(D'^{\,\pm1})=D'^{\,\pm1}\otimes D'^{\,\pm1},\\
\Delta(w_i)=w_i\ot w_i, \qquad \Delta(w_i')=w_i'\ot w_i',\\
\Delta(e_i)=e_i\ot 1+w_i\ot e_i, \qquad \Delta(f_i)=f_i\ot w_i'+1\ot
f_i,\\
\vep(e_i)=\vep(f_i)=0,\quad \vep(\gamma^{\pm\frac{1}2})
=\vep(\gamma'^{\,\pm\frac{1}2})=\vep(D^{\pm1})=\vep(D'^{\,\pm1})=\vep(w_i)=\vep(w_i')=1,
\\
S(\gamma^{\pm\frac{1}2})=\gamma^{\mp\frac{1}2},\qquad
S(\gamma'^{\pm\frac{1}2})=\gamma'^{\mp\frac{1}2},\qquad
S(D^{\pm1})=D^{\mp1},\qquad S(D'^{\,\pm1})=D'^{\,\mp1},\\
S(e_i)=-w_i^{-1}e_i,\qquad S(f_i)=-f_i\,w_i'^{-1},\qquad
S(w_i)=w_i^{-1}, \qquad S(w_i')=w_i'^{-1}.
\end{gather*}
\end{defi}

\noindent {\bf 2.2} \ In what follows, we give the skew-pairing and
the Drinfel'd double structure.
\begin{defi}
A bilinear form $\langle\, , \rangle:$ ${\frak{B}}\times
{\frak{A}}\longrightarrow \mathbb{K}$ is called a skew-dual pairing
of two Hopf algebras ${\frak{A}}$ and ${\frak{B}}$ $($see
$[\textrm{KS}, 8.2.1]$ $)$, if it satisfies
\begin{gather*}
\langle
b,\,
1_{{\frak{A}}}\rangle=\varepsilon_{{\frak{B}}}(b),\qquad\qquad
\langle 1_{{\frak{B}}},\, a\rangle=\varepsilon_{{\frak{A}}}(a),\\
\langle b,\, a_1a_2\rangle=\langle\Delta^{{\rm op}}_{{\frak{B}}}(b),
a_1\otimes a_2\rangle,\qquad \langle b_1b_2, a\rangle=\langle
b_1\otimes b_2, \Delta_{{\frak{A}}}(a)\rangle,
\end{gather*}
for all $a,\, a_1,\, a_2\in{{\frak{A}}}$ and $b,\, b_1,\,
b_2\in{{\frak{B}}}$, where $\varepsilon_{{\frak{A}}}$,
$\varepsilon_{{\frak{B}}}$ denote the counites of ${\frak{A}}$,
${\frak{B}}$, respectively, and $\Delta_{{\frak{A}}}$,
$\Delta_{{\frak{B}}}$ are the respective coproducts.
\end{defi}
\begin{defi} For any two Hopf algebras ${\frak{A}}$ and
${\frak{B}}$ skew-paired by $\langle\, , \rangle$, there exists a
Drinfel'd quantum double ${\cal D}({\frak{A}}, {\frak{B}})$ which is
a Hopf algebra whose underlying coalgebra is
${\frak{A}}\otimes{\frak{B}}$ with the tensor product coalgebra
structure, whose algebra structure is defined by
$$
(a\otimes b)(a'\otimes b')=\sum\langle S_{\frak{B}}(b_{(1)}),
a'_{(1)} \rangle\langle b_{(3)}, a'_{(3)}\rangle aa'_{(2)}\otimes
b_{(2)}b', $$
 for $a, a' \in {\frak{A}}$ and $b, b' \in {\frak{B}}$, and whose
 antipode $S$ is given by
$$
S(a\otimes b)=(1\otimes S_{\frak{B}}(b))(S_{\frak{A}}(a) \otimes 1).
$$
\end{defi}
Let $\widehat{\cal B}$ (resp. $\widehat{\cal B'}$) denote the Hopf
(Borel-type) subalgebra of $U_{r, s}(\widehat{\frak {sl}_n})$
generated by $e_j,~ \omega_j^{ \pm1}$, $\gamma^{\pm\frac{1}2}$,
$D^{\pm1}$ (resp. $f_j,\, \omega_j'^{\, \pm1}, \,
\gamma'^{\pm\frac{1}2},\, D'^{\,\pm1}$) with $j\in I_0$.

\begin{prop} There exists a unique skew-dual pairing
$\langle \,, \rangle:\widehat{\cal B'} \times \widehat{\cal
B}\longrightarrow\mathbb{K}$ of the Hopf subalgebras
$\widehat{\cal B}$ and $\widehat{\cal B'}$ such that:
\begin{gather*}
\langle f_i,\, e_j\rangle=\delta_{ij}\frac{1}{s-r},  \qquad\qquad (i,\,j\in I_0)\tag{1}\\
 \langle \omega'_i,\,\omega_j\rangle
=\begin{cases}r^{(\varepsilon_j,\,
\alpha_i)}\,s^{(\varepsilon_{j+1},\, \alpha_i)}, \quad\ \; (i\in
I_0,\, j\in I)\cr r^{-(\varepsilon_{i+1},\,
\alpha_0)}\,s^{(\varepsilon_1,\, \alpha_i)}, \quad (i\in
I_0,\,j=0)\end{cases}\tag{2}\\
\langle\omega_i'^{\,\pm1},\,
\omega_j^{-1}\rangle=\langle\omega_i'^{\, \pm1},\,
\omega_j\rangle^{-1} =\langle\omega'_i,\, \omega_j\rangle^{\mp1},
\quad (i,\,j\in I_0) \tag{3}\\
\langle\gamma'^{\,\pm\frac{1}{2}},\, \gamma\rangle=\langle\gamma',\,
\gamma^{\pm\frac{1}{2}}\rangle=\langle\gamma',\,
\gamma\rangle^{\pm\frac{1}2}=1,\tag{4}\\
\langle D',\, D\rangle^{\pm1}=\langle D'^{\,\pm1},\,
D\rangle=\langle D',\, D^{\pm1}\rangle=1,\tag{5}\\
\langle \gamma'^{\,\pm\frac{1}2},\, \om_i^{\pm1}\rangle=1=\langle
\om_i'^{\,\pm1},\, \gamma^{\pm\frac{1}2}\rangle,\quad (i\in I_0)
\tag{6}\\
\langle D'^{\pm1},\om_i\rangle=\langle D',
\om_i^{\pm1}\rangle=s^{\mp\delta_{0i}}, \quad \langle\om_i'^{\pm1},
D\rangle=\langle\om_i',
D^{\pm1}\rangle=r^{\pm\delta_{0i}}, (i\in I_0)\tag{7}\\
 \langle D',
\gamma^{\pm\frac{1}2}\rangle=\langle D'^{\,\pm1},
\gamma^{\frac{1}2}\rangle=s^{\mp\frac{1}2},\qquad\langle
\gamma'^{\,\pm\frac{1}2}, D\rangle=\langle \gamma'^{\frac{1}2},
D^{\pm1}\rangle=r^{\pm\frac{1}2},\tag{8}
\end{gather*}
and all others pairs of generators are $0$. Moreover, we have
$\langle S(b'), S(b)\rangle=\langle b', b\rangle$ for $b'\in
\widehat{\cal B'}, ~b\in\widehat{\cal B}.$
\end{prop}
\begin{proof} \ The uniqueness assertion is clear, as any
skew-dual pairing of bialgebras is determined by the values on the
generators. We proceed to prove the existence of the pairing.

The pairing defined on generators as (1)---(8) may be extended to a
bilinear form on $\widehat{\cal B'} \times \widehat{\cal B}$ in a
way such that the defining properties in Definition 2.2 hold. We
will verify that the relations in $\widehat{\cal B}$ and
$\widehat{\cal B'}$ are preserved, ensuring that the form is
well-defined and is a skew-dual pairing of $\widehat{\cal B}$ and
$\widehat{\cal B'}$.

First, it is straightforward to check that the bilinear form
preserves all the relations among the $\om_i^{\pm1}$,
$\gamma^{\pm\frac{1}2}$, $D^{\pm1}$ in $\widehat{\cal B}$ and the
${\om_i'}^{\pm1}$, $\gamma'^{\,\pm\frac{1}2}$, $D'^{\,\pm1}$ in
$\widehat{\cal B'}$. Next, we observe that the identities hold: for
$\,i,\, j\in I$,
$$
(\varepsilon_j, \alpha_i)=-(\varepsilon_{i+1}, \alpha_j), \qquad
(\varepsilon_j,\alpha_0)=-(\varepsilon_1, \alpha_j), \leqno(2.1)
$$
which ensure the compatibility of the form defined above with the
relations of $(\textrm{A2})$ and $(\textrm{A3})$ in $\widehat{\cal
B}$ or $\widehat{\cal B'}$ respectively. This fact is easily checked
by definition (see (1)---(8)).  So we are left to verify that the
form preserves the $(r,s)$-Serre relations in $\widehat{\cal B}$ and
$\widehat{\cal B'}$.

For $1 \leq i <n$, $(r,s)$-Serre relations in $\widehat{\cal B}$
and $\widehat{\cal B'}$ have been checked in [BW1]. Here we need
only to verify the relations involving index $i=0$ in
$\widehat{\cal B}$ and $\widehat{\cal B'}$. It suffices to
consider the following case (the remaining case is similar)
$$\langle X, e_0^2e_{n-1}-(r^{-1}+s^{-1})e_0e_{n-1}e_0+(rs)^{-1}e_{n-1}e_0^2\rangle,$$
where $X$ is any word in the generators of  $\widehat{\cal B'}$. By
definition, this equals
\begin{equation*}
\begin{split}
\langle \Delta^{(2)}(X),\,& \, e_0\otimes e_0\otimes
e_{n-1}\\
&-(r^{-1}+s^{-1})e_0\otimes e_{n-1}\otimes e_0+(r
s)^{-1}e_{n-1}\otimes e_0 \otimes e_0\rangle,
\end{split}\tag{2.2}
\end{equation*}
where $\Delta$ stands for $\Delta_{\cal B'}^{{\rm op}}$. In order
for any one of these terms to be nonzero, $X$ must involve exactly
two $f_0$ factors, one $f_{n-1}$ factor, and arbitrarily many
$\omega_j'^{\pm1}$ $(j\in I_0)$, $\gamma'^{\,\pm\frac{1}2}$, or
$D'^{\,\pm1}$ factors. For simplicity, we first consider three key
cases:

(i) \ If $X=f_0^2f_{n-1}$, then $\Delta^{(2)}(X)$ is equal to
\begin{eqnarray*}
&(\omega'_0\otimes \omega'_0 \otimes f_0+\omega'_0\otimes f_0
\otimes 1+ f_0 \otimes 1 \otimes 1)^2(\omega'_{n-1}\otimes
\omega'_{n-1}\otimes f_{n-1}\\
&\quad+\omega'_{n-1}\otimes f_{n-1}\otimes 1+f_{n-1}\otimes 1\otimes
1 ).
\end{eqnarray*}
 The relevant terms of
$\Delta^{(2)}(X)$ are
\begin{gather*}
f_0\omega'_0\omega'_{n-1}\otimes f_0\omega'_{n-1}\otimes f_{n-1}+
\omega'_0f_0\omega'_{n-1} \otimes
f_0\omega'_{n-1}\otimes f_{n-1}\\
+ f_0\omega'_0\omega'_{n-1}\otimes \omega'_0f_{n-1}\otimes f_0+
\omega'_0 f_0 \omega'_{n-1}\otimes \omega'_0 f_{n-1}\otimes f_0\\
+ \omega_0'^2f_{n-1}\otimes f_0 \omega'_0\otimes f_0+
\omega_0'^2f_{n-1}\otimes \omega'_0f_0 \otimes f_0.
\end{gather*}
Therefore, $(2.2)$ becomes
\begin{equation*}
\begin{split}
\langle f_0\,&\omega'_0\omega'_{n-1}, e_0 \rangle\langle
f_0\omega'_{n-1}, e_0 \rangle\langle f_{n-1}, e_{n-1}
\rangle\\
&\quad+\langle \omega'_0f_0\omega'_{n-1} , e_0\rangle\langle
 f_0\omega'_{n-1}, e_0\rangle \langle f_{n-1}, e_{n-1}\rangle \\
&\quad-(r^{-1}{+}s^{-1})\bigl(\langle f_0\omega'_0\omega'_{n-1},
e_0\rangle
 \langle \omega'_0f_{n-1}, e_{n-1}\rangle\langle
f_0, e_0 \rangle\\
&\quad+\langle \omega'_0 f_0 \omega'_{n-1}, e_0 \rangle\langle
\omega'_0f_{n-1}, e_{n-1} \rangle\langle f_0,
e_0\rangle\bigr)\\
&\quad +(rs)^{-1} \bigl(\langle {\omega}_0'^2f_{n-1}, e_{n-1}
\rangle\langle f_0 \omega_0', e_0\rangle\langle f_0, e_0
\rangle\\
&\quad+\langle
 {\omega}_0'^2f_{n-1}, e_{n-1}\rangle\langle \omega_0'f_0, e_0 \rangle\langle f_0,e_0\rangle\bigr)\\
&=\frac{1}{(s-r)^3}\bigl\{1+\langle\omega'_0,
\omega_0\rangle-(r^{-1}{+}s^{-1})\bigl(
 \langle \omega'_0, \omega_{n-1}\rangle+
 \langle \omega'_0, \omega_0 \rangle \langle \omega'_0, \omega_{n-1}\rangle\bigr)\\
&\quad+
 (rs)^{-1}\bigl(\langle \omega'_0, \omega_{n-1}\rangle^2+
 \langle\omega'_0,\omega_{n-1}\rangle^2\langle\omega'_0,
 \omega_0\rangle\bigr)\bigr\}\\
&=\frac{1}{(s-r)^3}\bigl\{1+rs^{-1}-(r^{-1}{+}s^{-1})(s+rs^{-1}s)+(rs)^{-1}(s^2+
 s^2rs^{-1})\bigr\}\\
&=0.
\end{split}
\end{equation*}

(ii) \ When $X=f_0f_{n-1}f_0$, it is easy to get the relevant terms
of $\Delta^{(2)}(X)$:
\begin{gather*}
\omega'_0\omega'_{n-1}f_0\otimes f_0\omega'_{n-1}\otimes f_{n-1}+
f_0\omega'_{n-1}\omega'_0\otimes \omega'_{n-1}f_0\otimes
f_{n-1}\\
+\omega'_0\omega'_{n-1}f_0\otimes \omega'_0f_{n-1}\otimes f_0 +f_0
\omega'_{n-1}\omega'_0\otimes f_{n-1}\omega'_0 \otimes f_0\\
+\omega'_0f_{n-1}\omega'_0 \otimes \omega'_0f_0\otimes f_0+
\omega'_0 f_{n-1} \omega'_0\otimes f_0\omega'_0 \otimes f_0.
\end{gather*}
Thus, $(2.2)$ becomes
\begin{equation*}
\begin{split}
&\frac{1}{(s{-}r)^3}\bigl\{\langle\omega'_0,
\omega_0\rangle\langle\omega'_{n-1}, \omega_0\rangle
 +\langle\omega'_{n-1}, \omega_0\rangle\\
&\quad-(r^{-1}{+}s^{-1}) \bigl(\langle \omega'_0, \omega_0\rangle
 \langle\omega'_{n-1}, \omega_0\rangle\langle\omega'_0, \omega_{n-1}\rangle
 +1\bigr)\\
&\quad +(r s)^{-1}\bigl(\langle \omega'_0,
\omega_{n-1}\rangle\langle\omega'_0, \omega_0\rangle+
 \langle\omega'_0,\omega_{n-1}\rangle\bigr)\bigr\}\\
 &=\frac{1}{(s{-}r)^3}\bigl\{rs^{-1}\cdot r^{-1}{+}r^{-1}{-}(r^{-1}{+}s^{-1})(rs^{-1}\cdot
 r^{-1}s{+}1){+}
 (r s)^{-1}(s\cdot rs^{-1}{+}s)\bigr\}\\
& =0.
\end{split}
\end{equation*}

(iii) \ If $X=f_{n-1}f_0^2$, one can similarly get that $(2.2)$
vanishes.

\medskip
Finally, if $X$ is any word involving exactly two $f_0$ factors, one
$f_{n-1}$ factor, and arbitrarily many factors $\omega_j'^{\, \pm1}$
$(j\in I_0$),  $\gamma'^{\,\pm\frac{1}2}$ and $D'^{\,\pm1}$, then
$(2.2)$ will just be a scalar multiple of one of the quantities we
have already calculated, and then will be 0.

Analogous calculations show that the relations in $\widehat{\cal
B'}$ are preserved.
\end{proof}

\begin{theorem}
${\cal D}({\widehat{\cal B}}, {\widehat{\cal B'}})$ is isomorphic to
$U_{r,s}(\widehat{\frak {sl}_n})$ as Hopf algebras.
\end{theorem}
\begin{proof} We denote the image $e_i\otimes 1$ of
$e_i$ in ${\cal D}({\widehat{\cal B}}, {\widehat{\cal B'}})$ by
${\widehat{e_i}}$ and similarly for $\omega_i^{\pm1},
\gamma^{\pm\frac{1}2}, D^{\pm1}$, denote the image $1\otimes f_i$ of
$f_i$ in ${\cal D}({\widehat{\cal B}}, {\widehat{\cal B'}})$ by
${\widehat{f_i}}$ and similarly for $\omega_i'^{\,\pm1},
\gamma'^{\,\pm\frac{1}2}, D'^{\,\pm1}$. Define a map $\varphi: {\cal
D}({\widehat{\cal B}}, {\widehat{\cal B'}})\longrightarrow
U_{r,s}(\widehat{\frak{sl}_n})$ by
\begin{gather*}
\varphi({\widehat{e_i}})=e_i,\quad\varphi({\widehat{f_i}})=f_i,\quad
\varphi(\widehat{\omega}_i^{\pm1})=\omega_i^{\pm1},\quad\varphi(\widehat{\omega}_i'^{
\pm1})=\omega_i'^{\pm1},\\
\varphi({\widehat{\gamma}^{\pm\frac{1}{2}}})=\gamma^{\pm\frac{1}{2}},
\quad\varphi({\widehat{\gamma'}^{\pm\frac{1}{2}}})={\gamma'}^{\pm\frac{1}{2}},
\quad
\varphi(\widehat{D}^{\pm1})=D^{\pm1},\quad\varphi(\widehat{D'}^{
\pm1})={D'}^{\pm1}.
\end{gather*}
The remaining argument is analogous to that of [BGH1, Theorem 2.5].
\end{proof}

\begin{remark} \ (1) Up to now, we have completely solved
the compatibility problem on the defining relations of our
two-parameter quantum affine algebra $U_{r,s}(\widehat{\frak
{sl}_n})$ $(n>2)$. This is done in two steps:
 the proof of Theorem 2.5 indicates that the cross
relations between $\widehat{\cal B}$ and $\widehat{\cal B'}$ are
half of the relations $(\textrm{A1})$---$(\textrm{A4})$, and the
proof of Proposition 2.4 shows the remaining relations,  including
the remaining half of relations $(\textrm{A1})$---$(\textrm{A4})$
and the $(r,s)$-Serre relations $(\textrm{A5})$---$(\textrm{A7})$.

\smallskip
(2) When $r=s^{-1}=q$, the Hopf algebra $U_{q,
q^{-1}}(\widehat{\frak {sl}_n})$ modulo the Hopf ideal generated by
the set $\{\,\om_i'-\om_i^{-1}$ $(i\in I_0)$,
$\gamma'^{\,\frac{1}2}-\gamma^{-\frac{1}2}$, $D'-D^{-1}\,\}$ is the
usual quantum affine algebra $U_q(\widehat{\frak {sl}_n})$ of
Drinfel'd-Jimbo type.
\end{remark}

Let $U^0=\mathbb
K[\om_0^{\pm1},\cdots,\om_n^{\pm1},{\om_0'}^{\pm1},\cdots,{\om_n'}^{\pm1}]$,
$U_0=\mathbb K[\om_0^{\pm1},\cdots,\om_n^{\pm1}]$, and $U_0'=\mathbb
K[{\om_0'}^{\pm1},\cdots,{\om_n'}^{\pm1}]$ denote the Laurent
polynomial subalgebras of $U_{r,s}(\widehat{\frak {sl}_n})$,
$\widehat{\mathcal B}$,
 and $\widehat{\mathcal B'}$ respectively. Clearly, $U^0=U_0U_0'=U_0'U_0$.
  Furthermore, let us
denote by $U_{r,s}(\widehat{\frak n})$ $($resp.
$U_{r,s}(\widehat{\frak n}^-)$\,$)$ the subalgebra of
$\widehat{\mathcal B}$ $($resp. $\widehat{\mathcal B'})$ generated
by $e_i$ $($resp. $f_i$$)$ for all $i\in I_0$. Thus, by
definition, we have $\widehat{\mathcal B}=U_{r,s}(\widehat{\frak
n})\rtimes U_0$, and $\widehat{\mathcal B'}=U_0'\ltimes
U_{r,s}(\widehat{\frak n}^-)$, so that the double $\mathcal
D(\widehat{\mathcal B},\widehat{\mathcal B'})\cong
U_{r,s}(\widehat{\frak n})\ot U^0\ot U_{r,s}(\widehat{\frak
n}^-)$, as vector spaces. On the other hand, if we consider $\la\,
,\, \ra^-: \widehat{\mathcal B'}\times \widehat{\mathcal
B}\lra\mathbb K$ by $\la b', b\ra^-:=\la S(b'), b\ra$, the
convolution inverse of the skew-dual paring $\la\, ,\ra$ in
Proposition 2.4, the composition with the flip mapping $\sigma$
then gives rise to a new skew-dual paring
$\la\,|\,\ra:=\la\,,\ra^-\circ\sigma: \widehat{\mathcal
B}\times\widehat{\mathcal B'}\lra \mathbb K$, given by $\la
b|b'\ra=\la S(b'),b\ra$. As a byproduct of Theorem 2.5,  similar
to [BGH1, Coro. 2.6], we get the standard triangular decomposition
of $U_{r,s}(\widehat{\frak{sl}_n})$.

\begin{coro}
$U_{r,s}(\widehat{\frak{sl}_n})\cong U_{r,s}(\widehat{\frak n}^-)\ot
U^0\ot U_{r,s}(\widehat{\frak n})$, as vector spaces.\hfill\qed
\end{coro}

\begin{coro} \ For any $\zeta=\sum_{i=0}^n\zeta_i\alpha_i\in
Q$ (the root lattice of $\widehat{\frak{sl}_n}$), the defining
relations (A2) and (A3) in $U_{r,s}(\widehat{\frak {sl}_n})$ take
the form:
\begin{gather*}
\om_{\zeta}\,e_i\,\om_{\zeta}^{-1}=\la \om_i',\om_\zeta\ra\,e_i,
\qquad
\om_{\zeta}\,f_i\,\om_{\zeta}^{-1}=\la \om_i',\om_\zeta\ra^{-1}f_i,\\
{\om_{\zeta}'}\,e_i\,{\om_{\zeta}'}^{-1}=\la \om_\zeta',
\om_i\ra^{-1} e_i,\qquad\quad
\om_{\zeta}'\,f_i\,{\om_{\zeta}'}^{-1}=\la \om_\zeta',\om_i\ra\,f_i.
\end{gather*}
$U_{r,s}(\widehat{\frak n}^\pm)=\bigoplus_{\eta\in
Q^+}U_{r,s}^{\pm\eta}(\widehat{\frak n}^\pm)$ is then $Q^\pm$-graded
with
$$
U_{r,s}^{\eta}({\widehat{\frak n}}^\pm)=\left\{\,a\in
U_{r,s}(\widehat{\frak n}^\pm)\;\left|\;
\om_\zeta\,a\,\om_\zeta^{-1}=\la \om_\eta',\om_\zeta\ra\,a, \
\om_\zeta'\,a\,{\om_\zeta'}^{-1}=\la \om_\zeta',\om_\eta\ra^{-1}
\,a\,\right\}\right.,
$$
for $\eta\in Q^+\cup Q^-$.

\medskip
Furthermore, $U=\bigoplus_{\eta\in
Q}U_{r,s}^\eta(\widehat{\frak{sl}_n})$ is $Q$-graded with
\begin{equation*}
\begin{split}
U_{r,s}^\eta(\widehat{\frak{sl}_n})&=\left\{\,\left.\sum
F_\alpha\om_{\mu}'\om_\nu E_\beta\in U\; \right|\;
\om_\zeta\,(F_\alpha\om_{\mu}'\om_\nu E_\beta)\,\om_{\zeta}^{-1}=
\la \om'_{\beta-\alpha},\om_\zeta\ra\,F_\alpha\om_\mu'\om_\nu E_\beta,\right.\\
&\ \left.\om_\zeta'\,(F_\alpha\om_\mu'\om_\nu
E_\beta)\,{\om_{\zeta}'}^{-1}= \la
\om_\zeta',\om_{\beta-\alpha}\ra^{-1}\,F_\alpha\om_\mu'\om_\nu
E_\beta,\; \textrm{\it with } \; \beta-\alpha=\eta\right\},
\end{split}
\end{equation*}
where $F_\alpha$ $($resp. $E_\beta$$)$ runs over monomials
$f_{i_1}{\cdots} f_{i_l}$ $($resp. $e_{j_1}{\cdots} e_{j_m}$$)$ such
that $\alpha_{i_1}+{\cdots}+\alpha_{i_l}=\alpha$ $($resp.
$\alpha_{j_1}+{\cdots}+\alpha_{j_m}=\beta$$)$.\hfill\qed
\end{coro}

\begin{defi}
Let $\tau$ be the $\mathbb{Q}$-algebra anti-automorphism of
$U_{r,s}(\widehat{\frak {sl}_n})$ such that $\tau(r)=s$,
$\tau(s)=r$, $\tau(\la \om_i',\om_j\ra^{\pm1})=\la
\om_j',\om_i\ra^{\mp1}$, and
\begin{gather*}
\tau(e_i)=f_i, \quad \tau(f_i)=e_i, \quad \tau(\om_i)=\om_i',\quad
\tau(\om_i')=\om_i,\\
\tau(\gamma)=\gamma',\quad
\tau(\gamma')=\gamma,\quad\tau(D)=D',\quad \tau(D')=D.
\end{gather*}
Then ${\widehat{\cal B'}}=\tau({\widehat{\cal B}})$ with those
induced defining relations from ${\widehat{\cal B}}$, and those
cross relations in $(\textrm{A2})$---$(\textrm{A4})$ are
antisymmetric with respect to $\tau$.\hfill\qed
\end{defi}

\medskip
\section{Drinfel'd
Realization of $U_{r,s}(\widehat{\frak{sl}_n})$ and Quantum Affine
Lyndon Basis}
\medskip

\noindent {\bf 3.1} \ For the two-parameter quantum affine algebra
$U_{r,s}(\widehat{\frak {sl}_n})$ $(n>2)$ we defined in Section 2,
we give the following definition of its Drinfel'd realization. In
the two-parameter case, the defining relations $(\textrm{D2})$,
$(\textrm{D6})$, $(\textrm{D7})$ and $(\textrm{D8})$ below appear to
vary dramatically in comparison with the one-parameter cases (see
$(\textrm{d2})$, $(\textrm{d6})$, $(\textrm{d7})$ and
$(\textrm{d8})$ in Remark 3.3), where the compatibilities for the
whole system are based on some intrinsic considerations as indicated
in the sequel.

We briefly write $\la i,j\ra:=\la \om_i',\om_j\ra$.
\begin{defi} Let ${\cal U}_{r,s}(\widehat{\frak {sl}_n})$ $(n>2)$ be the unital
associative algebra over $\mathbb{K}$ generated by the elements
$x_i^{\pm}(k)$, $a_i(\ell)$, $\om_i^{\pm1}$, ${\om'_i}^{\pm1}$,
$\gamma^{\pm\frac{1}{2}}$, ${\gamma'}^{\,\pm\frac{1}2}$, $D^{\pm1}$,
$D'^{\,\pm1}$ $(i\in I$, $k,\,k' \in \mathbb{Z}$, $\ell,\,\ell' \in
\mathbb{Z}\backslash \{0\})$, subject to the following defining
relations:

\medskip
\noindent $(\textrm{D1})$ \  $\gamma^{\pm\frac{1}{2}}$,
$\gamma'^{\,\pm\frac{1}{2}}$ are central with $\gamma\gamma'=rs$,
$\omega_i\,\omega_i^{-1}=\omega_j'\,\omega_j'^{\,-1}=1=DD^{-1}=D'D'^{-1}$
$(i, j\in I)$, and
\begin{equation*}
\begin{split}
[\,\omega_i^{\pm 1},\omega_j^{\,\pm 1}\,]&=[\,\om_i^{\pm1},
D^{\pm1}\,]=[\,\om_j'^{\,\pm1}, D^{\pm1}\,] =[\,\om_i^{\pm1},
D'^{\pm1}\,]=0\\
&=[\,\omega_i^{\pm 1},\omega_j'^{\,\pm 1}\,]=[\,\om_j'^{\,\pm1},
D'^{\pm1}\,]=[D'^{\,\pm1}, D^{\pm1}]=[\,\omega_i'^{\pm
1},\omega_j'^{\,\pm 1}\,].
\end{split}
\end{equation*}

$$[\,a_i(\ell),a_j(\ell')\,]=\delta_{\ell+\ell',0}\frac{(rs)^{\frac{|\ell|}2}
(\langle i,i\rangle^{\frac{\ell a_{ij}}2}-\langle
i,i\rangle^{-\frac{\ell a_{ij}}2})}{|\ell|(r-s)}
\cdot\frac{\gamma^{|\ell|}-\gamma'^{|\ell|}}{r-s}.
\leqno(\textrm{D2})$$
$$[\,a_i(\ell),~\om_j^{{\pm }1}\,]=[\,\,a_i(\ell),~{\om'}_j^{\pm
1}\,]=0.\leqno(\textrm{D3})
$$
\begin{gather*}
D\,x_i^{\pm}(k)\,D^{-1}=r^k\, x_i^{\pm}(k), \qquad\ \
D'\,x_i^{\pm}(k)\,D'^{\,-1}=s^k\, x_i^{\pm}(k),\tag{\textrm{D4}}
\\
D\, a_i(\ell)\,D^{-1}=r^\ell\,a_i(\ell), \qquad\qquad D'\,
a_i(\ell)\,D'^{\,-1}=s^\ell\,a_i(\ell).
\end{gather*}
$$
\om_i\,x_j^{\pm}(k)\, \om_i^{-1} =  \langle j, i\rangle^{\pm 1}
x_j^{\pm}(k), \qquad \om'_i\,x_j^{\pm}(k)\, \om_i'^{\,-1} = \langle
i, j\rangle ^{\mp1}x_j^{\pm}(k).\leqno(\textrm{D5})
$$
\begin{gather*}
[\,a_i(\ell),x_j^{\pm}(k)\,]=\pm\frac{(rs)^{\frac{|\ell|}2}(\langle
i,i\rangle^{\frac{\ell a_{ij}}2}-\langle i,i\rangle^{-\frac{\ell
a_{ij}}2})}{\ell(r-s)}\gamma^{\pm\frac{\ell}2}x_j^{\pm}(\ell{+}k),
\quad \textit{for} \
\ell<0,\tag{\textrm{D$6_1$}} \\
[\,a_i(\ell),x_j^{\pm}(k)\,]=\pm\frac{(rs)^{\frac{|\ell|}2}(\langle
i,i\rangle^{\frac{\ell a_{ij}}2}-\langle i,i\rangle^{-\frac{\ell
a_{ij}}2})}{\ell(r-s)}\gamma'^{\pm\frac{\ell}2}x_j^{\pm}(\ell{+}k),
\quad \textit{for} \ \ell>0.\tag{\textrm{D$6_2$}}
\end{gather*}

\begin{equation*}
\begin{split}
x_i^{\pm}(k{+}1)\,x_j^{\pm}(k') - \la j,i\ra^{\pm1} &x_j^{\pm}(k')\,x_i^{\pm}(k{+}1)\\
=-\Bigl(\la j,i\ra\la
i,j\ra^{-1}\Bigr)^{\pm\frac1{2}}\,&\Bigl(x_j^{\pm}(k'{+}1)\,x_i^{\pm}(k)-\la
i,j\ra^{\pm1} x_i^{\pm}(k)\,x_j^{\pm}(k'{+}1)\Bigr).
\end{split}\tag{\textrm{D7}}
\end{equation*}

$$
[\,x_i^{+}(k),~x_j^-(k')\,]=\frac{\delta_{ij}}{r-s}\Big(\gamma'^{-k}\,{\gamma}^{-\frac{k+k'}{2}}\,
\om_i(k{+}k')-\gamma^{k'}\,\gamma'^{\frac{k+k'}{2}}\,\om'_i(k{+}k')\Big),\leqno(\textrm{D8})
$$
where $\om_i(m)$, $\om'_i(-m)~(m\in \mathbb{Z}_{\geq 0})$ with
$\om_i(0)=\om_i$ and $\om'_i(0)=\om_i'$ are defined by:
\begin{gather*}
\sum\limits_{m=0}^{\infty}\om_i(m) z^{-m}=\om_i \exp \Big(
(r{-}s)\sum\limits_{\ell=1}^{\infty}
 a_i(\ell)z^{-\ell}\Big); \\
\sum\limits_{m=0}^{\infty}\om'_i(-m) z^{m}=\om'_i \exp
\Big({-}(r{-}s)
\sum\limits_{\ell=1}^{\infty}a_i(-\ell)z^{\ell}\Big),
\end{gather*}
with $\om_i(-m)=0$ and $\om'_i(m)=0, \ \forall\;m>0$.

$$x_i^{\pm}(m)x_j^{\pm}(k)=x_j^{\pm}(k)x_i^{\pm}(m),
\qquad\ \textit{for} \quad a_{ij}=0,\leqno(\textrm{D$9_1$})$$
\begin{equation*}
\begin{split}
& Sym_{m_1,\, m_2}\Big(x_i^{\pm}(m_1)
x_i^{\pm}(m_2)x_j^{\pm}(k)-(r^{\pm1}+s^{\pm1})\,x_i^{\pm}(m_1)x_j^{\pm}(k)x_i^{\pm}(m_2)\\
&\ +(rs)^{\pm1} x_j^{\pm}(k)x_i^{\pm}(m_1)x_i^{\pm}(m_2)\Big)=0,
\quad\textit{for}\quad a_{ij}=-1, \  1\leq i<j<n,\\
\end{split} \tag{\textrm{D$9_2$}}
\end{equation*}
\begin{equation*}
\begin{split}
& Sym_{m_1,\, m_2}\Big(x_i^{\pm}(m_1)
x_i^{\pm}(m_2)x_j^{\pm}(k)-(r^{\mp1}+s^{\mp1})\,x_i^{\pm}(m_1)x_j^{\pm}(k)x_i^{\pm}(m_2)\\
&\ +(rs)^{\mp1} x_j^{\pm}(k)x_i^{\pm}(m_1)x_i^{\pm}(m_2)\Big)=0,
\quad\textit{for}\quad a_{ij}=-1, \ 1\leq j<i<n,
\end{split}\tag{\textrm{D$9_3$}}
\end{equation*}
$\textit{Sym}$ denotes symmetrization with respect to the indices
$(m_1, m_2)$.
\end{defi}

As one of crucial observations of the compatibilities of the
defining system above, we have
\begin{prop}
There exists the $\mathbb{Q}$-algebra antiautomorphism $\tau$ of
$\,{\cal U}_{r,s}(\widehat{\frak {sl}_n})$ $(n>2)$ such that
$\tau(r)=s$, $\tau(s)=r$,
$\tau(\la\om_i',\om_j\ra^{\pm1})=\la\om_j',\om_i\ra^{\mp1}$ and
\begin{gather*}
\tau(\om_i)=\om_i',\quad \tau(\om_i')=\om_i,\\
\tau(\gamma)=\gamma',\quad \tau(\gamma')=\gamma,\\
\tau(D)=D',\quad \tau(D')=D,\\
\tau(a_i(\ell))=a_i(-\ell),\\
\tau(x_i^{\pm}(m))=x_i^{\mp}(-m), \\
\tau(\om_i(m))=\om_i'(-m),
\quad\tau(\om_i'(-m))=\om_i(m),
\end{gather*}
and $\tau$ preserves each defining relation $($\hbox{{\rm D}n}$)$ in
Definition 3.1 for $n=1,\cdots,9$.\hfill\qed
\end{prop}

\begin{remark} \ (1) Note that the defining relations
$(\textrm{D1})$---$(\textrm{D5})$, $(\textrm{D7})$, $(\textrm{D8})$,
and $(\textrm{D$9_1$})$---$(\textrm{D$9_3$})$ are self-compatible
each under the $\mathbb{Q}$-algebra antiautomorphism $\tau$, while
the couple of the defining relations
($(\textrm{D$6_1$})$,$(\textrm{D$6_2$})$) is compatible with each
other with respect to $\tau$. Using such a $\tau$, it is sufficient
to consider the compatibility for half of the relations, e.g., those
relations involving in $+$-parts for $x_i^{\pm}(m)$, or in positive
$\ell$'s for $a_i(\ell)$ (for instance, see $(\textrm{D$6_2$})$).

\smallskip
(2) The constraint condition $\gamma\gamma'=rs$ in $(\textrm{D1})$
is required intrinsically by the compatibilities among
$(\textrm{D1})$, $(\textrm{D3})$, $(\textrm{D5})$, $(\textrm{D6})$,
$(\textrm{D7})$ \& $(\textrm{D8})$. For instance, by
$(\textrm{D7})$, we have $[\,x_i^-(0), x_j^-(1)\,]_{\langle
i,j\rangle}=(\langle j,i\rangle\langle
i,j\rangle)^{\frac{1}2}[\,x_i^-(1),x_j^-(0)\,]_{\langle
j,i\rangle^{-1}}$. Thus, using the property (3.5) in Definition 3.4
below and $(\textrm{D8})$ \& $(\textrm{D5})$, we get
\begin{equation*}
\begin{split}
[\,x_j^+(0),[\,x_i^-(0), x_j^-(1)\,]_{\langle i, j\rangle}\,]&=
(\langle j,i\rangle\langle
i,j\rangle)^{\frac{1}2}[\,x_i^-(1),[\,x_j^+(0),x_j^-(0)\,]\,]_{\langle
j,i\rangle^{-1}}\\
&=(\langle j,i\rangle\langle
i,j\rangle)^{\frac{1}2}\Bigl[\,x_i^-(1),\frac{\om_j-\om_j'}{r-s}\,\Bigr]_{\langle
j,i\rangle^{-1}}\\
&=\frac{(\langle j,i\rangle\langle i,j\rangle)^{\frac{1}2}-(\langle
j,i\rangle\langle i,j\rangle)^{-\frac{1}2}}{r-s}x_i^-(1)\om_j.
\end{split}
\end{equation*}
However, using (3.5), $(\textrm{D8})$, $(\textrm{D3})$,
$(\textrm{D5})$ \& $(\textrm{D$6_2$})$, we can follow another way to
expand $[\,x_j^+(0),[\,x_i^-(0), x_j^-(1)\,]_{\langle i,
j\rangle}\,]$ directly as
\begin{equation*}
\begin{split}
[\,x_j^+(0),[\,x_i^-(0), x_j^-(1)\,]_{\langle i, j\rangle}\,]&=
[\,x_i^-(0),[\,x_j^+(0),x_j^-(1)\,]\,]_{\langle
i,j\rangle}\\
&=\gamma^{-\frac1{2}}[\,x_i^-(0),a_j(1)\,]\,\om_j\\
&=(rs)^{\frac1{2}}(\gamma\gamma')^{-\frac1{2}}\frac{\langle
j,j\rangle^{\frac{a_{ji}}2}-\langle
j,j\rangle^{-\frac{a_{ji}}2}}{r-s}x_i^-(1)\om_j.
\end{split}
\end{equation*}
Therefore, we obtain that $\gamma\gamma'=rs$ and $\langle
i,j\rangle\langle j,i\rangle=\langle i,i\rangle^{a_{ij}}$, for any
$i,\, j\in I$.

\smallskip
(3) As a glimpse of the compatibility of $(\textrm{D2})$ with
$(\textrm{D$6_1$})$, $(\textrm{D$6_2$})$ and $(\textrm{D8})$, we
have the following: By $(\textrm{D8})$, we get
$a_i(1)=\om_i^{-1}\gamma^{\frac1{2}}\,[\,x_i^+(0),\,x_i^-(1)\,]$ and
$a_i(-1)=\om_i'^{-1}{\gamma'}^{\frac1{2}}\,[\,x_i^+(-1),\,x_i^-(0)\,]$.
Then using one of these expressions of $a_i(\pm 1)$ and using
$(\textrm{D$6_1$})$ (or $(\textrm{D$6_2$})$) and $(\textrm{D8})$
again, we may expand the Lie bracket $[a_i(1), a_j(-1)]$ in two
manners to get to the same formula as $(\textrm{D2})$. One is to
expand $a_i(1)$ first, and then to use $(\textrm{D$6_1$})$ \&
$(\textrm{D8})$ as follows
\begin{equation*}
\begin{split}
[\,a_i(1), a_j(-1)\,]&=\om_i^{-1}\gamma^{\frac{1}{2}}[\,[x_i^+(0),
x_i^-(1)], a_j(-1)\,]\\
&=\om_i^{-1}\gamma^{\frac{1}2}\Bigl([\,[x_i^+(0), a_j(-1)],
x_i^-(1)\,]+[\,x_i^+(0), [x_i^-(1), a_j(-1)]\,]\Bigr)\\
&=\om_i^{-1}\gamma^{\frac{1}2}\{-a_{ij}\}_{\la
i,i\ra}\Bigl(\gamma^{-\frac{1}{2}}[x_i^+(-1),
x_i^-(1)]-\gamma^{\frac{1}{2}}[x_i^+(0),x_i^-(0)]\Bigr)\\
&=\{-a_{ij}\}_{\la
i,i\ra}\om_i^{-1}\Bigl(\frac{\gamma'\om_i-\gamma\om_i'}{r-s}
-\gamma\frac{\om_i-\om_i'}{r-s}\Bigr)\\
&=\{-a_{ij}\}_{\la i,i\ra}\frac{\gamma'-\gamma}{r-s}=\{a_{ij}\}_{\la
i,i\ra}\frac{\gamma-\gamma'}{r-s}.
\end{split}
\end{equation*}
where $\{\ell a_{ij}\}_{\la i,i\ra}:=\frac{(rs)^{\frac{|\ell|}2}(\la
i,i\ra^{\frac{\ell a_{ij}}2}-\la i,i\ra^{-\frac{\ell
a_{ij}}2})}{\ell(r-s)}=\{\ell a_{ji}\}_{\la j,j\ra}$,
$\{-a_{ij}\}_{\la i,i\ra}=-\{a_{ij}\}_{\la i,i\ra}$.
 Expanding $a_j(-1)$ instead and using $(\textrm{D$6_2$})$ \&
$(\textrm{D8})$, we get the same result. More compatibilities will
be clearer in the proof of the Drinfel'd isomorphism theorem.

\smallskip
(4) Another observation is the following: When $r=s^{-1}=q$, the
algebra ${\cal U}_{q, q^{-1}}(\widehat{\frak {sl}_n})$ modulo the
ideal generated by the set $\{\,\om_i'-\om_i^{-1}$ $(i\in I)$,
$\gamma'^{\,\frac{1}2}-\gamma^{-\frac{1}2}$, $D'-D^{-1}\,\}$ is
exactly the usual Drinfel'd realization ${\cal U}_q(\widehat{\frak
{sl}_n})$ defined below (cf. [B2]).

The unital associative algebra ${\cal U}_q(\widehat{\frak {sl}_n})$
over $\mathbb{Q}(q)$ is generated by the elements $x_i^{\pm}(k)$,
$a_i(\ell)$, $\om_i^{\pm1}$, $\gamma^{\pm\frac{1}{2}}$, $D^{\pm1}$,
$(i\in I$, $k \in \mathbb{Z}$, $\ell \in \mathbb{Z}\backslash
\{0\})$ subject to the following defining relations:

\medskip
\noindent $(\textrm{d1})$ \  $\gamma^{\pm\frac{1}{2}}$ are
central, $\omega_i\,\omega_i^{-1}=1=DD^{-1}$ $(i\in I)$, and for
$i,\,j\in I$, one has
$$[\,\omega_i^{\pm 1},\omega_j^{\,\pm 1}\,]=[\,\om_i^{\pm1},
D^{\pm1}\,] =0.
$$
$$[\,a_i(\ell),~a_j(\ell')\,]=\delta_{\ell+\ell',\,0}\frac{[\,\ell\,a_{ij}\,]}{\ell}\cdot
\frac{\gamma^{\ell}-\gamma^{-\ell}}{q-q^{-1}}, \qquad
\left([\,n\,]=\frac{q^n-q^{-n}}{q-q^{-1}}\right).
\leqno(\textrm{d2})
$$
$$[\,a_i(\ell),~\om_j^{{\pm }1}\,]=0.\leqno(\textrm{d3})
$$
$$
D\,x_i^{\pm}(k)\,D^{-1}=q^k\, x_i^{\pm}(k), \qquad D\,
a_i(\ell)\,D^{-1}=q^\ell\,a_i(\ell).\leqno(\textrm{d4})
$$
$$
\om_i\,x_j^{\pm}(k)\, \om_i^{-1} = q^{\pm a_{ij}}
x_j^{\pm}(k).\leqno(\textrm{d5})
$$
$$
[\,a_i(\ell),x_j^{\pm}(k)\,]=\pm\frac{[\,\ell\,
a_{ij}\,]}{\ell}\,\gamma^{\mp\frac{|\ell|}{2}}\,x_j^{\pm}(\ell{+}k).
 \leqno(\textrm{d$6$})
$$
\begin{equation*}
\begin{split}
x_i^{\pm}(k+1)x_j^{\pm}(k') &- q^{\pm a_{ij}} x_j^{\pm}(k')x_i^{\pm}(k{+}1)\\
&=q^{\pm a_{ij}} x_i^{\pm}(k)x_j^{\pm}(k'{+}1)-
x_j^{\pm}(k'{+}1)x_i^{\pm}(k).
\end{split}\tag{\textrm{d7}}
\end{equation*}
$$
[\,x_i^{+}(k),~x_j^-(k')\,]=\frac{\delta_{ij}}{q-q^{-1}}\Big(\gamma^{\frac{k-k'}{2}}\,
\om_i(k{+}k')-\gamma^{\frac{k'-k}{2}}\,\om^{-1}_i(k{+}k')\Big),\leqno(\textrm{d8})
$$
where $\om_i(m)$ and $\om^{-1}_i(-m)~(m\in \mathbb{Z}_{\geq 0})$
with $\om_i(0)=\om_i$ and $\om^{-1}_i(0)=\om_i^{-1}$ are defined by:
\begin{gather*}\sum\limits_{m=0}^{\infty}\om_i(m) z^{-m}=\om_i \exp \Big(
(q{-}q^{-1})\sum\limits_{\ell=1}^{\infty}
 a_i(\ell)z^{-\ell}\Big),\quad \bigl(\om_i(-m)=0, \ \forall\;m>0\bigr); \\
\sum\limits_{m=0}^{\infty}\om^{-1}_i(-m) z^{m}=\om^{-1}_i \exp
\Big({-}(q{-}q^{-1})
\sum\limits_{\ell=1}^{\infty}a_i(-\ell)z^{\ell}\Big), \quad
\bigl(\om^{-1}_i(m)=0, \ \forall\;m>0\bigr).
\end{gather*}
$$x_i^{\pm}(m)x_j^{\pm}(k)=x_j^{\pm}(k)x_i^{\pm}(m),
\qquad\ \textit{for} \quad a_{ij}=0,\leqno(\textrm{d$9_1$})$$
\begin{equation*}
\begin{split}
& Sym_{m_1,\, m_2}\Big(x_i^{\pm}(m_1)
x_i^{\pm}(m_2)x_j^{\pm}(k)-(q^{\pm1}+q^{\mp1})\,x_i^{\pm}(m_1)x_j^{\pm}(k)x_i^{\pm}(m_2)\\
&\ +x_j^{\pm}(k)x_i^{\pm}(m_1)x_i^{\pm}(m_2)\Big)=0,
\quad\textit{for}\quad a_{ij}=-1, \  1\leq i<j<n,\\
\end{split} \tag{\textrm{d$9_2$}}
\end{equation*}
\begin{equation*}
\begin{split}
& Sym_{m_1,\, m_2}\Big(x_i^{\pm}(m_1)
x_i^{\pm}(m_2)x_j^{\pm}(k)-(q^{\mp1}+q^{\pm1})\,x_i^{\pm}(m_1)x_j^{\pm}(k)x_i^{\pm}(m_2)\\
&\ + x_j^{\pm}(k)x_i^{\pm}(m_1)x_i^{\pm}(m_2)\Big)=0,
\quad\textit{for}\quad a_{ij}=-1, \ 1\leq j<i<n.
\end{split}\tag{\textrm{d$9_3$}}
\end{equation*}
\end{remark}

\noindent {\bf 3.2} \ Before putting forward the Drinfel'd
isomorphism theorem, that is, showing that the
$\mathbb{Q}(r,s)$-algebra ${\cal U}_{r,s}(\widehat{\frak {sl}_n})$
$(n>2)$ in Definition 3.1 is exactly the Drinfel'd realization of
the two-parameter quantum affine algebra $ U_{r,s}(\widehat{\frak
{sl}_n})$ $(n>2)$ defined in Definition 2.1, we need to make some
preliminaries on Lyndon words, and to adapt a definition of quantum
Lie bracket borrowed from [J2] to give our definition about
``affine" quantum Lie bracket (see Definition 3.6) which enables us
to derive an interesting description on the quantum {\it affine}
Lyndon basis in the quantum affine cases for the first time.

Note that the (affine) quantum Lie bracket possesses some advantages
in calculations such as less related to degrees of elements (see the
properties (3.3) \& (3.4) below). This generalized quantum Lie
bracket, like the one used in the usual construction of the quantum
Lyndon basis (for definition, see [R2]), is consistent with the
process when adding the bracketing on those corresponding Lyndon
words. This is crucial to the quantum calculations we develop later
on.

\begin{defi} $([\hbox{\rm J}2])$ \ The quantum Lie bracket $\,[\,a_1, a_2,\cdots,
a_s\,]_{(q_1,\,q_2,\,\cdots,\, q_{s-1})}$  is  defined inductively
by
\begin{gather*} [\,a_1, a_2\,]_q=a_1a_2-q\,a_2a_1,\quad
\hbox{\it for \ $q\in \mathbb{K}\backslash \{0\}$},\\
[\,a_1, a_2, \cdots, a_s\,]_{(q_1,\,q_2,\cdots,q_{s-1})}=[\,a_1,
[\,a_2, \cdots, a_s\,]_{(q_1,\cdots,q_{s-2})}\,]_{q_{s-1}},
\quad\hbox{\it for \ $q_i\in \mathbb{K}\backslash \{0\}$}.
\end{gather*}

The following identities follow from the definition.
\begin{gather*}
[\,a, bc\,]_v=[\,a, b\,]_x\,c+x\,b\,[\,a, c\,]_{\frac{v}{x}}, \qquad x\neq 0,\tag{3.1}\\
[\,ab, c\,]_v=a\,[\,b, c\,]_x+x\,[\,a, c\,]_{\frac{v}{x}}\,b, \qquad x\neq 0,\tag{3.2}\\
[\,a,[\,b,c\,]_u\,]_v=[\,[\,a,b\,]_x,
c\,]_{\frac{uv}{x}}+x\,[\,b,[\,a,c\,]
_{\frac{v}{x}}\,]_{\frac{u}{x}},
\qquad x\neq 0, \tag{3.3}\\
[\,[\,a,b\,]_u,c\,]_v=[\,a,[\,b,c\,]_x\,]_{\frac{uv}{x}}+x\,[\,[\,a,c\,]
_{\frac{v}{x}},b\,]_{\frac{u}{x}}, \qquad x\neq 0.\tag{3.4}\\
\bigl[\,a, [\,b_1, \cdots, b_s\,]_{(v_1,\,\cdots,\,
v_{s-1})}\,\bigr] =\sum_i\bigl[\,b_1,\cdots,[\,a, b_i\,],
\cdots,b_s\,\bigr]_{(v_1,\,\cdots,\, v_{s-1})},\tag{3.5}\\
[\,a, a, b\,]_{(u,\,
v)}=[\,a,a,b\,]_{(v,\,u)}=a^2b-(u+v)\,aba+(uv)\,ba^2.\tag{3.6}
\end{gather*}
\end{defi}

\begin{defi} For the generators system of the algebra
$\mathcal{U}_{r,s}(\widehat{\frak{sl}_n})$, we define the $\dot
Q$-gradation (where $\dot Q$ is the root lattice of $\frak{sl}_n$)
as follows:
\begin{gather*}
\textrm{deg}(\om_i^{\pm1})=\textrm{deg}(\om_i'^{\pm1})=\textrm{deg}(\gamma^{\pm\frac1{2}})
=\textrm{deg}(\gamma'^{\pm\frac{1}2})=\textrm{deg}(D^{\pm1})
=\textrm{deg}(D'^{\pm1})=0,\\
\textrm{deg}(a_i(\pm\ell))=0, \qquad
\textrm{deg}(x_i^{\pm}(k))=\pm\alpha_i.
\end{gather*}
Hence, the defining relations (D1)---(D9) ensure that
$\mathcal{U}_{r,s}(\widehat{\frak{sl}_n})$ has a triangular
decomposition:
$$\mathcal{U}_{r,s}(\widehat{\frak{sl}_n})=
\mathcal{U}_{r,s}(\widetilde{\frak{n}}^-)\otimes\mathcal{U}_{r,s}^0(\widehat{\frak{sl}_n})
\otimes\mathcal{U}_{r,s}(\widetilde{\frak{n}}),$$ where
$\mathcal{U}_{r,s}(\widetilde{\frak{n}}^\pm)=\bigoplus_{\alpha\in\dot
Q^\pm}\mathcal{U}_{r,s}(\widetilde{\frak{n}}^\pm)_\alpha$ is
generated respectively by $x_i^\pm(k)$ ($i\in I$), and
$\mathcal{U}_{r,s}^0(\widehat{\frak{sl}_n})$ is the subalgebra
generated by $\om_i^{\pm1}$, $\om_i'^{\pm1}$,
$\gamma^{\pm\frac1{2}}$, $\gamma'^{\pm\frac1{2}}$, $D^{\pm1}$,
$D'^{\pm1}$ and $a_i(\pm\ell)$ for $i\in I$, $\ell\in \mathbb{N}$.
Namely, $\mathcal{U}_{r,s}^0(\widehat{\frak{sl}_n})$ is generated by
the toral subalgebra $\mathcal{U}_{r,s}(\widehat{\frak{sl}_n})^0$
and the quantum Heisenberg subalgebra $\mathcal
H_{r,s}(\widehat{\frak{sl}_n})$ generated by those quantum imaginary
root vectors $a_i(\pm\ell)$ ($i\in I$, $\ell\in \mathbb{N}$).
\end{defi}

\begin{defi}
For $\alpha,\,\beta\in\dot Q^+$ (a positive root lattice of
$\frak{sl}_n$),  $x_\alpha^\pm(k),\,x_\beta^\pm(k')$ $\in
\mathcal{U}_{r,s}(\widetilde{\frak n}^\pm)$, we define their
 {\it ``affine" quantum Lie bracket} as follows:
$$\bigl[\,x_\alpha^\pm(k),\,
x_\beta^\pm(k')\,\bigr]_{\la\om'_\alpha,\om_\beta\ra^{\mp1}}:=
x_\alpha^\pm(k)\,x_\beta^\pm(k')-\la\om'_\alpha,\om_\beta\ra^{\mp1}x_\beta^\pm(k')\,x_\alpha^\pm(k).\leqno(3.7)
$$
\end{defi}
By definition 3.6, the formula $(\textrm{D7})$ will take the
convenient form as
$$\left[\,x^{\pm}_i(k),\,
x_j^{\pm}(k'{+}1)\,\right]_{\langle
i,j\rangle^{\mp1}}=-\Bigl(\langle j,i\rangle\langle
i,j\rangle^{-1}\Bigr)^{\pm\frac1{2}} \left[\,x^{\pm}_j(k'),\,
x_i^{\pm}(k{+}1)\,\right]_{\langle j,i\rangle^{\mp1}}. \leqno(3.8)
$$

By (3.6), the $(r,s)$-Serre relations (D$9_2$) \& (D$9_3$) for
$m_1=m_2$ in the case of $a_{ij}=-1$ can be reformulated as:
\begin{equation*}
\begin{split}
\bigl[\,x_i^{\pm}(m), x_i^{\pm}(m),
x_j^{\pm}(k)\,\bigr]_{(r^{\pm1},\,s^{\pm1})}&=0,
\qquad\textit{for} \quad 1\leq i<j<n,\\
\bigl[\,x_i^{\pm}(m), x_i^{\pm}(m),
x_j^{\pm}(k)\,\bigr]_{(s^{\mp1},\,r^{\mp1})}&=0, \qquad
\textit{for}\quad 1\leq j<i<n.
\end{split}\tag{3.9}
\end{equation*}

\begin{remark} (1) \ For any nonsimple root $\alpha\, (\ne\alpha_i)$ ($i\in I$), the
meaning of notation $x_\alpha^+(k)$ (resp. $x_\alpha^-(k)$) in
Definition 3.6 has a bit ambiguity, as is well-known even for
quantum ``classical" root vectors $x_\alpha^+(0)$ which have
different linearly-independent choices. However, the combinatorial
approach to Lyndon words, together with the ``affine" quantum Lie
bracket, will give us a valid and specific choice for
$x_\alpha^+(k)$ which leads to a construction of quantum ``affine"
Lyndon basis for $\mathcal U_{r,s}(\widetilde{\frak n})$, on which
acting $\tau$ will yield a corresponding construction of quantum
``affine" Lyndon basis for $\mathcal U_{r,s}(\widetilde{\frak n}^-)$
(see Proposition 3.10 \& Theorem 3.11 below).

\medskip
(2) \ In fact, (3.8) describes a kind of consistent constraints of
quantum affine root vectors defined by some Lyndon words of
different levels (if say, $x_j^\pm(k)$ have {\it level} $k$) which
obeys the defining rule of Lyndon basis (see below) via Lyndon words
as in the classical types, since from (3.8), we get the {\it
level-shifting formula}
\begin{equation*}
\begin{split}
\bigl[\,x^{\pm}_i(k),\, x_j^{\pm}(k'{+}1)\,\bigr]_{\langle
i,j\rangle^{\mp1}} &=\la
i,i\ra^{\mp\frac{a_{ij}}2}\Bigl(\left[\,x^{\pm}_i(k{+}1),\,
x_j^{\pm}(k')\,\right]_{\langle
i,j\rangle^{\mp1}}\\
&\quad +\bigl(\la i,j\ra^{\mp1} {-}\la j,i\ra^{\pm1}
\bigr)\,x_j^{\pm}(k')\,x^{\pm}_i(k{+}1)\Bigr).
\end{split}\tag{3.10}
\end{equation*}
Based on this formula, we will see that it makes reasonable to give
the definition of quantum affine root vector $x_\alpha^\pm(k)$ as in
(3.14) \& (3.15) below such that the level $k$ completely
concentrates on the component of the lowest index, in the ordered
constituents of Lyndon basis. This will be clear from the proof of
Proposition 3.10.

\medskip
(3) \ Let $\mathcal U_{r,s}(\frak{n})$ denote the subalgebra of
$\mathcal U_{r,s}(\widetilde{\frak n})$, generated by $x_i^+(0)$
($i\in I$). By definition, it is clear that $\mathcal
U_{r,s}(\frak{n})\cong U_{r,s}(\frak n)$, the subalgebra of
$U_{r,s}(\frak{sl}_n)$ generated by $e_i$ ($i\in I$) (see [BGH1,
Remarks (2), p. 391]). Now let us recall the construction of a
Lyndon basis. The natural ordering $<$ in $I$ gives a total ordering
of the alphabet $A=\{x_1^+(0),\cdots,x_{n-1}^+(0)\}$. Let $A^*$ be
the set of all words in the alphabet $A$ (including the vacuum $1$)
and let $u<v$ denote that word $u$ is lexicographically smaller than
word $v$. Recall that a word $\ell \in A^*$ is a {\it Lyndon word}
if it is lexicographically smaller than all its proper right factors
(cf. [LR], [R2], [BH]). Let ${\mathbb{K}}[A^*]$ be the associative
algebra of $\mathbb{K}$-linear combinations of words in $A^*$ whose
product is juxtaposition, namely, a free ${\mathbb{K}}$-algebra. Let
$J$ be the $(r,s)$-Serre ideal of ${\mathbb{K}}[A^*]$ generated by
elements $\{(ad_\ell x_i^+(0))^{1-a_{ij}}(x_j^+(0))\mid 1\leq i\ne j
\leq n-1\}$. Clearly, $\mathcal
U_{r,s}({\mathfrak{n}})={\mathbb{K}}[A^*]/J$. Now given another
ordering $\preceq$ in $A^*$ with introducing a usual length function
$|\cdot|$ for each word $u\in A^*$. We say $u\preceq w $, if
$|u|<|w|$ or $|u|=|w|$ and $u\geq w$. Then we call a (Lyndon) word
to be {\it good } with respect to the $(r,s)$-Serre ideal $J$ if it
cannot be written as a sum of strictly smaller words modulo $J$ with
respect to the ordering $\preceq$. From [R2], the set of quantum Lie
brackets (or say, ${\mathbf q}$-bracketings) of all good Lyndon
words consists of a system of quantum root vectors of $\mathcal
U_{r,s}(\frak{n})$. More precisely, we have a construction for any
quantum root vector $x_\alpha^+(0)$ with $\alpha\in\dot\Delta^+$ (a
positive root system of $\frak{sl}_n$) in the following.\hfill\qed
\end{remark}

Take a corresponding ordering (compatible with the natural
ordering $<$ on $I$) of
$\dot\Delta^+=\{\alpha_{ij}:=\alpha_i+\alpha_{i+1}+\cdots+\alpha_{j-1}=\varepsilon_i-\varepsilon_j\mid
1\le i<j\le n\}$ with $\alpha_{i,i+1}=\alpha_i$ as follows (see
[H, p. 533]):
\begin{equation*}
\begin{split}
\alpha_{12},\, \alpha_{13},\, \alpha_{14},\cdots,\alpha_{1n},&\\
\alpha_{23},\, \alpha_{24},\,\cdots,\, \alpha_{2n},&\\
\vdots\quad &\\
 \alpha_{n-1,n}&
\end{split}\tag{3.11}
\end{equation*}
which is a {\it convex ordering} on
$\dot\Delta^+$ (for definition, see [R2, Section 6]). Hence, for
each $\alpha=\alpha_{ij}\in \dot\Delta^+$, by [R2], we can construct
the quantum root vector $x_\alpha^+(0)$ as a $(r,s)$-bracketing of a
good Lyndon word in the inductive fashion:
\begin{equation*}
\begin{split}
x_{\alpha_{ij}}^+(0):&=\bigl[\,x_{\alpha_{i,j-1}}^+(0), x_{j-1}^+(0)\,\bigr]_{\la \om_{\alpha_{i,j-1}}',\,\om_{j-1}\ra^{-1}}\\
&=\bigl[\,\cdots\bigl[\,x_i^+(0), x_{i+1}^+(0)\,\bigr]_{\la
i,i+1\ra^{-1}},\cdots,x_{j-1}^+(0)\,\bigr]_{\la
\om_{\alpha_{i,j-1}}',\,\om_{j-1}\ra^{-1}}\\
&=\bigl[\,\cdots\bigl[\,x_i^+(0),
x_{i+1}^+(0)\,\bigr]_r,\cdots,x_{j-1}^+(0)\,\bigr]_r.
\end{split}\tag{3.12}
\end{equation*}

Applying $\tau$ to (3.12), we can obtain the definition of quantum
root vector $x_{\alpha_{ij}}^-(0)$ as below:
$$x_{\alpha_{ij}}^-(0)=\tau\bigl(x_{\alpha_{ij}}^+(0)\bigr)=\bigl[\,x_{j-1}^-(0),\cdots,
\bigl[\,x_{i+1}^-(0), x_i^-(0)\,\bigr]_s\cdots\,\bigr]_s.
\leqno(3.13)$$

\begin{theorem}
$(\textrm{\rm i})$ \ The set
$$\left\{\,\left.x_{\alpha_{n{-}1,n}}^+(0)^{\ell_{n{-}1,n}}\cdots
x_{\alpha_{23}}^+(0)^{\ell_{23}}x_{\alpha_{1n}}^+(0)^{\ell_{1n}}
\cdots
x_{\alpha_{13}}^+(0)^{\ell_{13}}x_{\alpha_{12}}^+(0)^{\ell_{12}}\,\right|\,
\ell_{ij}\ge0\,\right\}$$ is a Lyndon basis of $\mathcal
U_{r,s}(\frak{n})$.

\smallskip
$(\textrm{\rm ii})$ \ The set
$$\left\{\,\left.x_{\alpha_{12}}^-(0)^{\ell_{12}}x_{\alpha_{13}}^-(0)^{\ell_{13}}\cdots
x_{\alpha_{1n}}^-(0)^{\ell_{1n}}x_{\alpha_{23}}^-(0)^{\ell_{23}}
\cdots x_{\alpha_{n-1,n}}^-(0)^{\ell_{n-1,n}}\,\right|\,
\ell_{ij}\ge0\,\right\}$$ is a Lyndon basis of $\mathcal
U_{r,s}(\frak{n^-})$.\hfill\qed
\end{theorem}

\begin{defi}
For $\alpha_{ij}\in\dot\Delta^+$, we define the {\it quantum affine
root vectors} $x_{\alpha_{ij}}^\pm(k)$ of nontrivial {\it level }$k$
by
\begin{gather*}
x_{\alpha_{ij}}^+(k):=\bigl[\,\cdots\bigl[\,x_i^+(k),
x_{i+1}^+(0)\,\bigr]_r,\cdots,x_{j-1}^+(0)\,\bigr]_r,\tag{3.14}\\
x_{\alpha_{ij}}^-(k):=\bigl[\,x_{j-1}^-(0),\cdots,
\bigl[\,x_{i+1}^-(0), x_i^-(k)\,\bigr]_s\cdots\,\bigr]_s, \tag{3.15}
\end{gather*}
where $\tau\bigl(x_{\alpha_{ij}}^\pm(\pm
k)\bigr)=x_{\alpha_{ij}}^\mp(\mp k)$.
\end{defi}

For each fixed $\alpha\in \dot Q^+$, let us denote by $\mathcal
U_{r,s}^{(k)}(\widetilde{\frak n})_\alpha$ the subspace of $\mathcal
U_{r,s}(\widetilde{\frak n})_\alpha$, consisting of elements of
level $k$. Hence, $\mathcal U_{r,s}(\widetilde{\frak
n})_\alpha=\bigoplus_{k\in\mathbb Z}\mathcal
U_{r,s}^{(k)}(\widetilde{\frak n})_\alpha$. When
$\alpha=\alpha_i\in\dot \Delta^+$ is a simple root, by definition,
$\dim\mathcal U_{r,s}^{(k)}(\widetilde{\frak n})_{\alpha_i}=1$ for
any level $k$. However, for any nonsimple root $\alpha\ne\alpha_i$
$(i\in I$), $\dim\mathcal U_{r,s}^{(k)}(\widetilde{\frak
n})_{\alpha}=\infty$ for any level $k$. In this case, given a
positive root $\alpha=\alpha_{ij}\in\dot \Delta^+$, we call a tuple
$(\beta_{j_1},\cdots,\beta_{j_\nu})$ ($\nu\ge1$) to be a partition
of root $\alpha_{ij}$ if $\beta_{j_1}<\cdots<\beta_{j_\nu}$ in the
ordering given in (3.11) such that
$\beta_{j_1}+\cdots+\beta_{j_\nu}=\alpha_{ij}$. If $\nu>1$, we say
this partition to be {\it proper}. Denote by ${\frak
P}^\circ(\alpha)$ the set of all proper partitions of root $\alpha$.
Obviously, we have $\mathcal U_{r,s}^{(k_\nu)}(\widetilde{\frak
n})_{\beta_{j_\nu}}\cdots \mathcal U_{r,s}^{(k_1)}(\widetilde{\frak
n})_{\beta_{j_1}}\subseteq \mathcal U_{r,s}^{(k)}(\widetilde{\frak
n})_{\alpha}$ if $k_1+\cdots+k_\nu=k$. Now we write
$$
\Omega_\alpha^{(k)}(\widetilde{\frak
n}):=\sum_{(\beta_{j_1},\cdots,\beta_{j_\nu})\in{\frak
P}^\circ(\alpha)\atop k_1+\cdots k_\nu=k}\mathcal
U_{r,s}^{(k_\nu)}(\widetilde{\frak n})_{\beta_{j_\nu}}\cdots
\mathcal U_{r,s}^{(k_1)}(\widetilde{\frak
n})_{\beta_{j_1}}\subseteq \mathcal U_{r,s}^{(k)}(\widetilde{\frak
n})_{\alpha}
$$
for the subspace of $\mathcal U_{r,s}^{(k)}(\widetilde{\frak
n})_{\alpha}$ spanned by basis elements' products of level $k$
from those proper partitions pertaining to $\alpha$. Using the
$\mathbb{Q}$-antiautomorphism $\tau$ on
$\Omega_\alpha^{(-k)}(\widetilde{\frak n})$, we get
$$\Omega_\alpha^{(k)}(\widetilde{\frak n}^-):=\tau\bigl(\Omega_\alpha^{(-k)}(\widetilde{\frak
n})\bigr).$$ Then we have the following description on $\mathcal
U_{r,s}^{(k)}(\widetilde{\frak n}^\pm)_{\alpha}$ for
$\alpha\in\dot\Delta^+$, whose proof shows that Definition 3.9
makes sense.
\begin{prop} For $1\le i<j\le n$ and $\alpha_{ij}\in\dot\Delta^+$
$($a positive root system of $\frak{sl}_n\,)$, we have

$(\textrm{\rm i})$ \ \;$\mathcal U_{r,s}^{(k)}(\widetilde{\frak
n})_{\alpha_{ij}}=\mathbb K
x_{\alpha_{ij}}^+(k)\bigoplus\Omega_{\alpha_{ij}}^{(k)}(\widetilde{\frak
n})$,

$(\textrm{\rm ii})$ \ $\mathcal U_{r,s}^{(k)}(\widetilde{\frak
n}^-)_{\alpha_{ij}}=\mathbb K
x_{\alpha_{ij}}^-(k)\bigoplus\Omega_{\alpha_{ij}}^{(k)}(\widetilde{\frak
n}^-)$.
\end{prop}
\begin{proof} \ (i) \
We will use an induction on rank $n$, where $n\ge 2$. Assume that
$i<j$ and $k'>0$, then by (3.10), we have
\begin{equation*}
\begin{split}
\bigl[\,x^+_i(k),\, x_j^+(k')\,\bigr]_{\la i,j\ra^{-1}} &=\la
i,i\ra^{-\frac{a_{ij}}2}\Bigl(\left[\,x^+_i(k{+}1),\,
x_j^+(k'{-}1)\,\right]_{\la i,j\ra^{-1}}\\
&\quad +\bigl(\la i,j\ra^{-1} {-}\la j,i\ra
\bigr)\,x_j^+(k'{-}1)\,x^+_i(k{+}1)\Bigr).
\end{split}\tag{3.16}
\end{equation*}
\begin{equation*}
\begin{split}
\bigl[\,x^+_i(k),\, x_j^+(-k')\,\bigr]_{\la i,j\ra^{-1}} &=\la
i,i\ra^{\frac{a_{ij}}2}\left[\,x^+_i(k{-}1),\,
x_j^+(-k'{+}1)\,\right]_{\la i,j\ra^{-1}}\\
&\quad +\bigl(\la j,i\ra-\la
i,j\ra^{-1}\bigr)\,x_j^+(-k')\,x^+_i(k).
\end{split}\tag{3.17}
\end{equation*}

When $n=2$, for any $k'\in\mathbb N$, repeatedly using (3.16) \&
(3.17), we get
\begin{equation*}
\begin{split}
\bigl[\,x^+_1(k),\, x_2^+(k')\,\bigr]_r &=\la
1,1\ra^{\frac{k'}2}\,x^+_{\alpha_{13}}(k{+}k')\\&\quad+\sum_{t=1}^{k'}\la
1,1\ra^{\frac{k'{-}t{+}1}2}(r{-}s)\,x_2^+(t{-}1)\,x_1^+(k{+}k'{-}t{+}1)\\
&\equiv \la 1,1\ra^{\frac{k'}2}\,x^+_{\alpha_{13}}(k{+}k')\quad
\mod\;
\Omega_{\alpha_{13}}^{(k+k')}(\widetilde{\frak n}),\\
 \bigl[\,x^+_1(k),\,
x_2^+(-k')\,\bigr]_r &=\la
1,1\ra^{-\frac{k'}2}\,x^+_{\alpha_{13}}(k{-}k')\\&\quad+\sum_{t=1}^{k'}\la
1,1\ra^{-\frac{k'{-}t}2}(s{-}r)\,x_2^+(-t)\,x_1^+(k{-}k'{+}t)\\&\equiv
\la 1,1\ra^{-\frac{k'}2}\,x^+_{\alpha_{13}}(k{-}k')\quad
\mod\;\Omega_{\alpha_{13}}^{(k-k')}(\widetilde{\frak n}),
\end{split}
\end{equation*}
which means that in both cases, we have
$$\bigl[\,x^+_1(k),\,
x_2^+(k')\,\bigr]_r \equiv \la
1,1\ra^{\frac{k'}2}\,x^+_{\alpha_{13}}(k{+}k')\quad \mod\;
\Omega_{\alpha_{13}}^{(k+k')}(\widetilde{\frak n}), \quad
\textit{for any } \ k'\in\mathbb Z.$$ Therefore, in rank $2$ case,
any elements (except for $x_2^+(k')x_1^+(k)$) of degree
$\alpha_{13}$ generated by $x_1^+(k)$ and $x_2^+(k')$ are of form:
$\bigl[\,x_1^+(k), x_2^+(k')\,\bigr]_a$ for any $a\in \mathbb K$,
however,
\begin{equation*}
\begin{split}
\bigl[\,x_1^+(k), x_2^+(k')\,\bigr]_a&=\bigl[\,x_1^+(k),
x_2^+(k')\,\bigr]_r+(r-a)\,x_2^+(k')\,x_1^+(k)\\
&\equiv (rs^{-1})^{\frac{k'}2}\,x^+_{\alpha_{13}}(k{+}k')\quad
\mod\; \Omega_{\alpha_{13}}^{(k+k')}(\widetilde{\frak n}).
\end{split}\tag{3.18}
\end{equation*}
This fact shows that
$$\mathcal U_{r,s}^{(k)}(\widetilde{\frak
n})_{\alpha_{13}}=\mathbb K
x_{\alpha_{13}}^+(k)\bigoplus\Omega_{\alpha_{13}}^{(k)}(\widetilde{\frak
n})$$
as vector spaces. Dually, we also have $\mathcal
U_{r,s}^{(k)}(\widetilde{\frak n}^-)_{\alpha_{13}}=\mathbb K
x_{\alpha_{13}}^-(k)\bigoplus\Omega_{\alpha_{13}}^{(k)}(\widetilde{\frak
n}^-)$ as vector spaces.

Now we assume that we have proved the results for rank $<n$, that
is, for those $\alpha_{ij}$ with $1\le i<j<n$. For rank $n$ case,
owing to the ordering given in (3.11), we are left to prove the
remaining cases: $\mathcal U_{r,s}^{(k)}(\widetilde{\frak
n}^\pm)_{\alpha_{in}}$ with $1\le i<j=n$.

In view of the same observation as (3.18), we need only to
consider the following elements of degree $\alpha_{in}$ and level
$k+k'$ generated by $x_{\alpha_{i,n-1}}^+(k)$ and $x_{n-1}^+(k')$
for $1\le i<n$: $\bigl[\,x_{\alpha_{i,n-1}}^+(k),
x_{n-1}^+(k')\,\bigr]_{\la \om'_{
\alpha_{i,n-1}},\,\om_{n-1}\ra^{-1}}=\bigl[\,x_{\alpha_{i,n-1}}^+(k),
x_{n-1}^+(k')\,\bigr]_r$. By definition (see (3.14)) and using
(3.4), (3.5) \& (3.1), we have
\begin{equation*}
\begin{split}
\bigl[\,x_{\alpha_{i,n-1}}^+&(k),
x_{n-1}^+(k')\,\bigr]_r=\bigl[\,\bigl[\,x_{\alpha_{i,n-2}}^+(k),x_{n-2}^+(0)\,\bigr]_r,
x_{n-1}^+(k')\,\bigr]_r\quad\textrm{(using (3.4))}\\
&=\bigl[\,x_{\alpha_{i,n-2}}^+(k),\bigl[\,x_{n-2}^+(0),
x_{n-1}^+(k')\,\bigr]_r\,\bigr]_r\\
&\quad+r\bigl[\,\underbrace{\bigl[\,x_{\alpha_{i,n-2}}^+(k),x_{n-1}^+(k')\,\bigr]},x_{n-2}^+(0)
\,\bigr]\\&\hskip2cm\textrm{($2$nd term $=0$ by (3.5) \&
(D9$_1$))}\\
&=\bigl[\,x_{\alpha_{i,n-2}}^+(k),\underbrace{\bigl[\,x_{n-2}^+(0),
x_{n-1}^+(k')\,\bigr]_r}\,\bigr]_r\quad\textrm{(using (3.18): rank $2$ case)}\\
&=(rs^{-1})^{\frac{k'}2}\bigl[\,\underbrace{\bigl[\,x_{\alpha_{i,n-2}}^+(k),x_{n-2}^+(k')\,\bigr]_r},
x_{n-1}^+(0)\,\bigr]_r\\
&\hskip2cm\textrm{(using the inductive hypothesis})\\
&\quad+\sum_{t}*_t\,(r-s)\,\bigl[\,x_{\alpha_{i,n-2}}^+(k),x_{n-1}^+(t)\,x_{n-2}^+(k'{-}t)\,\bigr]_r\quad\textrm{(using
(3.1))}\\
&\equiv
(rs^{-1})^{\frac{k'(n{-}1{-}i)}2}\underbrace{\bigl[\,x_{\alpha_{i,n-1}}^+(k{+}k'),x_{n-1}^+(0)\,\bigr]_r}
\ \;\mod\,\bigl[\,\Omega_{\alpha_{i,n-1}}^{(k{+}k')}(\widetilde{\frak{n}}),x_{n-1}^+(0)\,\bigr]_r\\
&\hskip3.5cm(\textrm{by definition})\\
 &\quad+\sum_{t}*_t\,(r-s)\,x_{n-1}^+(t)\,\underbrace{\bigl[\,x_{\alpha_{i,n-2}}^+(k),x_{n-2}^+(k'{-}t)\,\bigr]_r}\\
 &\hskip3.5cm(\textrm{using the inductive hypothesis})\\
&\equiv
(rs^{-1})^{\frac{k'(n{-}1{-}i)}2}x_{\alpha_{in}}^+(k{+}k')\quad\mod\,\Omega_{\alpha_{in}}^{(k{+}k')}(\widetilde{\frak{n}})\\
&\quad+\sum_{t}*_t'\,(r-s)\,x_{n-1}^+(t)\,x_{\alpha_{i,n-1}}^+(k{+}k'{-}t)\quad \mod\,x_{n-1}^+(t)\,\Omega_{\alpha_{i,n-1}}^{(k{+}k'{-}t)}(\widetilde{\frak{n}})\\
&\equiv
(rs^{-1})^{\frac{k'(n{-}1{-}i)}2}x_{\alpha_{in}}^+(k{+}k')\quad\mod\,\Omega_{\alpha_{in}}^{(k{+}k')}(\widetilde{\frak{n}}),
\end{split}
\end{equation*}
where in the $1$st ``$\equiv$", we used the following fact:
\begin{equation*}
\begin{split}
\bigl[\,x_{\alpha_{i,n-2}}^+(k),x_{n-1}^+(t)\,x_{n-2}^+(k'{-}t)\,\bigr]_r
&=x_{n-1}^+(t)\,\bigl[\,x_{\alpha_{i,n-2}}^+(k),x_{n-2}^+(k'{-}t)\,\bigr]_r\\
&\quad+\underbrace{\bigl[\,x_{\alpha_{i,n-2}}^+(k),x_{n-1}^+(t)\,\bigr]}\,x_{n-2}^+(k'{-}t)\\
&\hskip1cm(\textrm{$2$nd term $=0$ by (3.5) \& (D9$_1$)})
\\
&=x_{n-1}^+(t)\,\bigl[\,x_{\alpha_{i,n-2}}^+(k),x_{n-2}^+(k'{-}t)\,\bigr]_r;
\end{split}
\end{equation*}
while in the $2$nd ``$\equiv$", we used the facts:
\begin{gather*}
\bigl[\,\Omega_{\alpha_{i,n-1}}^{(k{+}k')}(\widetilde{\frak{n}}),x_{n-1}^+(0)\,\bigr]_r\subseteq
\Omega_{\alpha_{in}}^{(k{+}k')}(\widetilde{\frak{n}}),\\
x_{n-1}^+(t)\,\Omega_{\alpha_{i,n-1}}^{(k{+}k'{-}t)}(\widetilde{\frak{n}})\subseteq
\Omega_{\alpha_{in}}^{(k{+}k')}(\widetilde{\frak{n}}).
\end{gather*}
The latter is clear, due to the definition of
$\Omega_{\alpha_{in}}^{(k{+}k')}(\widetilde{\frak{n}})$. As for
the first inclusion, we have the following argument provided that
we notice the basis elements' constituents of
$\Omega_{\alpha_{i,n-1}}^{(k{+}k')}(\widetilde{\frak{n}})$.
Indeed, for any basis element
$$x^+_{\alpha_{\ell_\nu,n-1}}(k_\nu)\,x_{\alpha_{\ell_{\nu-1},\ell_{\nu}}}^+(k_{\nu-1})\cdots
x_{\alpha_{i,\ell_1}}^+(k_1)\in\Omega_{\alpha_{i,n-1}}^{(k{+}k')}(\widetilde{\frak{n}})
$$ of level $k+k'$ pertaining to a
partition of $\alpha_{i,n-1}$, using (3.2), we have
\begin{equation*}
\begin{split}
\bigl[\,x^+_{\alpha_{\ell_\nu,n-1}}&(k_\nu)\,x_{\alpha_{\ell_{\nu-1},\ell_{\nu}}}^+(k_{\nu-1})\cdots
x_{\alpha_{i,\ell_1}}^+(k_1),x_{n-1}^+(0)\,\bigr]_r\\
&=\underbrace{\bigl[\,x_{\alpha_{\ell_\nu,n-1}}^+(k_\nu),x_{n-1}^+(0)\,\bigr]_r}\,x_{\alpha_{\ell_{\nu-1},\ell_\nu}}^+(k_{\nu-1})\cdots
x_{\alpha_{i,\ell_1}}^+(k_1)\ \;(\textrm{by definition})\\
&\quad+x_{\alpha_{\ell_\nu,n-1}}^+(k_\nu)\underbrace{\bigl[\,x_{\alpha_{\ell_{\nu-1},\ell_\nu}}^+(k_{\nu-1})\cdots
x_{\alpha_{i,\ell_1}}^+(k_1),x_{n-1}^+(0)\,\bigr]}\\
&\quad(\textrm{$2$nd term $=0$ by (3.5) \& (D9$_1$) since $\ell_1<\cdots<\ell_\nu<n{-}1$})\\
&=x_{\alpha_{\ell_\nu,n}}^+(k_\nu)\,x_{\alpha_{\ell_{\nu-1},\ell_\nu}}^+(k_{\nu-1})\cdots
x_{\alpha_{i,\ell_1}}^+(k_1)\\
&\in\Omega_{\alpha_{in}}^{(k{+}k')}(\widetilde{\frak{n}}),
\end{split}
\end{equation*}
here $k_1+\cdots+k_\nu=k+k'$.

Up to now, we have finished the proof of (i). Using $\tau$ to (i),
we can get the second statement (ii).
\end{proof}

The argument above (in fact used the so-called quantum calculations)
implies the important conclusions about the quantum affine Lyndon
basis we present below.
\begin{theorem}
$(\textrm{\rm i})$ \ The set
$$\left\{\left.\left(\prod_{i\in\Bbb Z}x_{\alpha_{n{-}1,n}}^+(i)^{\ell^{(i)}_{n{-}1,n}}\right)\cdots
\left(\prod_{i\in\Bbb
Z}x_{\alpha_{1n}}^+(i)^{\ell^{(i)}_{1n}}\right) \cdots
\left(\prod_{i\in\Bbb
Z}x_{\alpha_{12}}^+(i)^{\ell^{(i)}_{12}}\right)\,\right|\,
\ell^{(i)}_{st}\ge0\,\right\}$$ is an ``affine" Lyndon basis of
$\mathcal U_{r,s}(\widetilde{\frak{n}})$, where each index set
$I_{\alpha_{st}}=\{i\in\Bbb Z\mid \ell^{(i)}_{st}\ne 0\}$ is finite.

\smallskip
$(\textrm{\rm ii})$ \ The set
$$\left\{\,\left.\left(\prod_{i\in\Bbb Z}x_{\alpha_{12}}^-(i)^{\ell^{(i)}_{12}}\right)\cdots
\left(\prod_{i\in\Bbb
Z}x_{\alpha_{1n}}^-(i)^{\ell^{(i)}_{1n}}\right) \cdots
\left(\prod_{i\in\Bbb
Z}x_{\alpha_{n-1,n}}^-(i)^{\ell^{(i)}_{n-1,n}}\right)\,\right|\,
\ell^{(i)}_{st}\ge0\,\right\}$$ is an ``affine" Lyndon basis of
$\mathcal U_{r,s}(\frak{\widetilde{n}^-})$, where each index set
$I_{\alpha_{st}}=\{i\in\Bbb Z\mid \ell^{(i)}_{st}\ne 0\}$ is
finite.\hfill\qed
\end{theorem}

\noindent {\bf 3.3} \ The following main theorem establishes the
Drinfel'd isomorphism between the two-parameter quantum affine
algebra $U_{r,s}(\widehat{\frak {sl}_n})$ (in Definition 2.1) and
the $(r,s)$-analogue of Drinfel'd quantum affinization of
$U_{r,s}(\frak{sl}_n)$ (in Definition 3.1), which affords the
two-parameter Drinfel'd realization of $U_{r,s}(\widehat{\frak
{sl}_n})$ as required.

\begin{theorem} $(${\bf Drinfel'd Isomorphism}$)$ \ For Lie algebra $\frak{sl}_n$ with $n>2$, let
$\theta=\alpha_{1n}$ be the maximal positive root. Then there exists
an algebra isomorphism $\Psi: U_{r,s}(\widehat{\frak {sl}_n})
\longrightarrow {\cal U}_{r,s}(\widehat{\frak {sl}_n})$ defined by:
for each $i\in I,$
\begin{eqnarray*}
\omega_i&\longmapsto& \om_i\\
\omega'_i&\longmapsto& \om'_i \\
\omega_0&\longmapsto& \gamma'^{-1}\, \om_{\theta}^{-1}\\
\omega'_0&\longmapsto& \gamma^{-1}\, \om_{\theta}'^{-1}\\
\end{eqnarray*}
\begin{eqnarray*}
\gamma^{\pm\frac{1}2}&\longmapsto& \gamma^{\pm\frac{1}2}\\
\gamma'^{\,\pm\frac{1}2}&\longmapsto&
\gamma'^{\,\pm\frac{1}2}\\
D^{\pm1}&\longmapsto& D^{\pm1}\\
D'^{\,\pm1}&\longmapsto& D'^{\,\pm1}\\
e_i&\longmapsto& x_i^+(0)\\
f_i&\longmapsto& x_i^-(0)\\
e_0&\longmapsto& x_{\alpha_{1n}}^-(1)
      \,\cdot(\gamma'^{-1}\,\om_{\theta}^{-1})=x^-_{\theta}(1)\cdot(\gamma'^{-1}\,\om_{\theta}^{-1})\\
f_0 &\longmapsto&
(\gamma^{-1}\,\om_{\theta}'^{-1})\cdot\,x^+_{\alpha_{1n}}(-1)=\tau\Bigl(x_{\alpha_{1n}}^-(1)
      \,\cdot(\gamma'^{-1}\,\om_{\theta}^{-1})\Bigr)
\end{eqnarray*}
where  $\om_{\theta}=\om_1\,\cdots\,
\om_{n-1},\,\om'_{\theta}=\om'_1\,\cdots\, \om'_{n-1}$.\hfill\qed
\end{theorem}

Since the Lusztig's symmetry of the braid group for the
two-parameter cases is no more available when the rank of $\frak g$
is bigger than $2$ (see [BGH1, Section 3]). This means that the
Beck's approach (using the extended braid group actions (see [B2])
to prove the Drinfel'd Isomorphism Theorem) is not yet valid for the
two-parameter cases here. Our treatment in the next section in fact
develops a valid and interesting algorithm on the quantum
calculations, which, as the reader has seen, is also a successful
application to the combinatorial approach to the quantum ``affine"
Lyndon basis (based on the Drinfel'd generators) we introduced
above. In some sense, our method also provides another new
combinatorial proof via the quantum ``affine" Lyndon basis even in
the one-parameter setting.

\medskip
\section{Proof of the Drinfel'd
Isomorphism Theorem}
\medskip

\noindent {\bf 4.1} \ Let $E_i,\,F_i$, $\om_i$, $\om_i'$ denote
the images of $e_i,\,f_i$, $\om_i$, $\om_i'$  ($i\in I_0$) in the
algebra ${\cal U}_{r,s}(\widehat{\frak {sl}_n})$ under the mapping
$\Psi$, respectively.

Denote by $\mathcal{U\,'}_{r,s}(\widehat{\frak{sl}_n})$ the
subalgebra of $\mathcal{U}_{r,s}(\widehat{\frak{sl}_n})$ generated
by $E_i,\,F_i,\,\om_i^{\pm1}$, $\om_i'^{\pm1}$ ($i\in I_0$),
$\gamma^{\pm\frac{1}2},\, \gamma'^{\,\pm\frac{1}2}$, $D^{\pm1}$
and $D'^{\pm1}$, that is,
$$
{\mathcal U\,'}_{r,s}(\widehat{\frak
{sl}_n}):=\left.\left\langle\, E_i,\, F_i,\, \om_i^{\pm1},\,
{\om}_i'^{\,\pm1}\;,\, \gamma^{\pm\frac{1}2},\,
\gamma'^{\,\pm\frac{1}2},\, D^{\pm1},\, D'^{\,\pm1}\; \right|
\;i\in I_0\;\right\rangle.
$$
Thereby, to prove the Drinfel'd Isomorphism Theorem (Theorem 3.12)
is equivalent to prove the following three Theorems:

\begin{theorem}
$\Psi:\,U_{r,s}(\widehat{\frak {sl}_n}) \longrightarrow {\mathcal
U\,'}_{r,s}(\widehat{\frak {sl}_n})$ is an epimorphism.
\end{theorem}

\begin{theorem}
${\mathcal U\,'}_{r,s}(\widehat{\frak {sl}_n})={\mathcal
U}_{r,s}(\widehat{\frak {sl}_n})$.
\end{theorem}

\begin{theorem}
$\Psi:\,U_{r,s}(\widehat{\frak {sl}_n}) \longrightarrow {\mathcal
U}_{r,s}(\widehat{\frak {sl}_n})$ is injective.
\end{theorem}

\medskip
\noindent {\bf 4.2} \ {\it Proof of Theorem 4.1.} \ We shall check
that the elements $E_i,~F_i,~\om_i,~\om'_i$ $(i\in I_0),\,
\gamma^{\pm\frac{1}2},\, \gamma'^{\,\pm\frac{1}2},\, D^{\pm1},
D'^{\pm1}$ satisfy the defining relations of
$(\textrm{A1})$--$(\textrm{A7})$ of $U_{r,s}(\widehat{\frak
{sl}_n})$. First of all, the defining relations of ${\cal
U}_{r,s}(\widehat{\frak {sl}_n})$ imply that
$E_i,\,F_i,\,\om_i,\,\om'_i \; (i\in I)$ generate a subalgebra
$\mathcal U_{r,s}(\frak{sl}_n)$ of ${\cal U}_{r,s}(\widehat{\frak
{sl}_n})$, which is isomorphic to $U_{r,s}(\frak{sl}_n)$. So we are
left to check the relations involving the index $i=0$.

\smallskip
Obviously, the relations of $(\textrm{A1})$ hold, according to the
defining relations of ${\cal U}_{r,s}(\widehat{\frak {sl}_n})$.

\medskip
For $(\textrm{A2})$: \ we just check the following three relations
involving $i=0$, the remaining relations in $(\textrm{A2})$ are
parallel to check. Using (\textrm{D}4), we get
$$
DE_0D^{-1}=D\,x^-_{\theta}(1)
\,D^{-1}\cdot(\gamma'^{-1}\om_{\theta}^{-1}) =r\,E_0.
$$
For $0\le j<n$, noting that $\langle\om_\theta'^{-1},
\om_j\rangle=\langle\gamma^{-1}\om_\theta'^{-1},
\om_j\rangle=\langle \om_0', \om_j\rangle$ (by Proposition 2.4), we
have
\begin{equation*}
\begin{split}
\om_j\,E_0\,\om_j^{-1}&=\om_j\,x^-_{\theta}(1) \,(\gamma'^{-1} \om_{\theta}^{-1})\,\om_j^{-1}\\
&=\langle \om_{n-1}', \om_j\rangle^{-1}\cdots\langle \om'_1,
\om_j\rangle^{-1}E_0=\langle\om'_0, \om_j\rangle E_0.
\end{split}
\end{equation*}
For $0\le i<n,$
$$\om_0\, E_i\, \om_0^{-1}=(\gamma'^{-1}\om_\theta^{-1})\,E_i\,(\gamma'\om_\theta)
=\om_\theta^{-1}\,E_i\,\om_\theta,$$ when $i\neq0$, since
$\langle\om_i', \om_\theta\rangle^{-1}=\langle\om_i',\om_0\rangle$
(by Proposition 2.4), we obtain
$$\om_\theta^{-1}\,E_i\,\om_\theta=\om_\theta^{-1}\,x_i^+(0)\,\om_\theta
=\langle \om'_i, \om_\theta\rangle^{-1}x_i^+(0)=\langle \om_i',
\om_0\rangle E_i;$$ and when $i=0$, since
$\langle\om_\theta',\om_\theta^{-1}\rangle^{-1}=\langle\om_0',\om_0\rangle$
(by Proposition 2.4), we have
$$
\om_\theta^{-1}\,E_0\,\om_\theta=\om_\theta^{-1}\,
x^-_{\theta}(1)\,(\gamma'^{-1}\om_{\theta}^{-1})\,\om_\theta
 =\langle\om_0',\om_0\rangle E_0.
$$

\smallskip
Similarly, one can verify the relations in $(\textrm{A3})$.

\medskip For $(\textrm{A4})$: \ first of all, when $i\neq 0$, we see
that
\begin{equation*}
\begin{split}
[\,E_0,F_i\,]
&=\bigl[\,x^-_{\theta}(1)\cdot(\gamma'^{-1}\om_{\theta}^{-1}),\,
x_i^-{(0)}\,\bigr]\\
&=-\bigl[\,x_i^-{(0)},\,x^-_{\theta}(1)\,\bigr]_{\langle \omega'_i,~
\omega_\theta \rangle}(\gamma'^{-1}\om_{\theta}^{-1}).
\end{split}
\end{equation*}

According to the corresponding cross relations held in
$U_{r,s}(\widehat{\frak{sl}_n})$, we claim the following crucial
Lemma, whose proof using the typical quantum calculations is
technical.

\begin{lemm} \ $\bigl[\,x_i^-{(0)},\,
 x^-_{\theta}(1)\,\bigr]_{\langle \omega_i',~
\omega_\theta\rangle}=0$,\quad for \ $i\in I$.
\end{lemm}
\begin{proof} \ (I) \ When $i=1$,
$\langle\om_1',\om_0\rangle=\langle\om_2',\om_1\rangle=s$, and $\la
\om_1',\om_\theta\ra=s^{-1}$. By (3.8) \& (3.9), we have
\begin{equation*}
\begin{split}
\bigl[\,x_1^-(0),\,x_{\alpha_{13}}^-(1)\,\bigr]_{s^{-1}}
&=\bigl[\,x_1^-(0),\,\underbrace{\bigl[\,x_2^-(0),\,x_1^-(1)\,\bigr]_{\langle2,1\rangle}}\,\bigr]_{s^{-1}}\qquad(\hbox{by (3.8)})\\
&=-(\langle1,2\rangle^{-1}\langle2,1\rangle)^{\frac1{2}}\,
\bigl[\,x_1^-(0),\,\bigl[\,x_1^-(0),\,x_2^-(1)\,\bigr]_{\langle1,2\rangle}\,\bigr]_{s^{-1}}\\
&=-(rs)^{\frac1{2}}\,
\bigl[\,x_1^-(0),\,x_1^-(0),\,x_2^-(1)\,\bigr]_{(r^{-1},s^{-1})}=0.\qquad(\hbox{by
(3.9)})
\end{split}
\end{equation*}
Hence, repeatedly using (3.3), we have
\begin{equation*}
\begin{split}
\bigl[\,x_1^-(0),\, x^-_{\alpha_{1n}}(1)\,\bigr]_{s^{-1}}
&=\bigl[\,\underbrace{\bigl[\,x_1^-(0),
x_{n-1}^-(0)\,\bigr]},\,x^-_{\alpha_{1,n-1}}(1)\,\bigr]\qquad(=0 \hbox{ by (\textrm{D9$_1$})})\\
&\quad+
\bigl[\,x_{n-1}^-(0),\, \bigl[\,x_1^-(0), x^-_{\alpha_{1,n-1}}(1)\,\bigr]_{s^{-1}}\,\bigr]_s\qquad(\hbox{by (3.3)})\\
 &=\bigl[\,x_{n-1}^-(0),\,\underbrace{
\bigl[\,x_1^-(0), x^-_{\alpha_{1,n-1}}(1)\,\bigr]_{s^{-1}}}\,\bigr]_s\\
&=\cdots\qquad(\hbox{inductively using (3.3) \& (\textrm{D9$_1$})})\\
&=\bigl[\,x_{n-1}^-(0),\,x_{n-2}^-(0),\cdots,\underbrace{
\bigl[\,x_1^-(0),\,x_{\alpha_{13}}^-(1)\,\bigr]_{s^{-1}}}\cdots\bigr]_{(s,\cdots,s)}\\
&=0.
\end{split}
\end{equation*}

(II) \ When $i=n-1$, $\langle\om_{n-1}',\om_0\rangle=r^{-1}$, that
is, $\la \om_{n-1}',\om_\theta\ra=r$. By (3.3), (3.9) \& (D9$_1$),
we have
\begin{equation*}
\begin{split}
\bigl[\,&x_{n-1}^-(0),\, x^-_{\alpha_{1n}}(1)\,\bigr]_r\qquad
 \hbox{(by definition)} \\
&=\bigl[\,x_{n-1}^-(0), \,\underbrace{\bigl[\,x_{n-1}^-(0),\,
\bigl[\,x_{n-2}^-(0),\,x^-_{\alpha_{1,n-2}}(1)\,\bigr]_s\bigr]_s}\,\bigr]_r \qquad(\hbox{using (3.3)})\\
 &=\underbrace{\bigl[\,x_{n-1}^-(0), \bigl[\,\bigl[\,x_{n-1}^-(0), x_{n-2}^-(0)\,\bigr]_{s},
 x^-_{\alpha_{1,n-2}}(1)\,\bigr]_s\,\bigr]_r}\qquad\hbox{(this term using (3.3))}\\
 &\quad +\,s\,\bigl[\,x_{n-1}^-(0), \bigl[\,x_{n-2}^-(0),\,
 \underbrace{\bigl[\,x_{n-1}^-(0), \,x^-_{\alpha_{1,n-2}}(1)\,\bigr]}\,\bigr]\,\bigr]_r\qquad(=0 \hbox{ by (3.5), (\textrm{D9$_1$})})\\
 &=\bigl[\,\underbrace{\bigl[\,x_{n-1}^-(0), x_{n-1}^-(0),
 x_{n-2}^-(0)\,\bigr]_{(s,r)}},\,x^-_{\alpha_{1,n-2}}(1)\,\bigr]_s\qquad(\hbox{this term$=0$ by (3.9))}\\
 &\quad+\,r\,
 \bigl[\,\bigl[\,x_{n-1}^-(0), x_{n-2}^-(0)\,\bigr]_s, \underbrace{\bigl[\,x_{n-1}^-(0),\,
 x^-_{\alpha_{1,n-2}}(1)\,\bigr]}\,\bigr]_{r^{-1}s}\quad(=0 \hbox{ by (3.5), (\textrm{D9$_1$})})\\
 &=0.
\end{split}
\end{equation*}

(III) \ When $1<i<n-1$, $\langle\om_i',\om_0\rangle=1$, that is,
$\langle\om_i',\om_\theta\rangle=1$. In order to derive the
required result, we first need to make two claims below:

\medskip
\noindent {\bf Claim $(A)$}: \
$\bigl[\,x_i^-(0),\,x^-_{\alpha_{1,i+1}}(1)\,\bigr]_{\la\om_i',\om_{\alpha_{1,i+1}}\ra}=\bigl[\,x_i^-(0),\,x^-_{\alpha_{1,i+1}}(1)\,\bigr]_{r}=0$,
\ \textit{for} \ $i\ge 2$.

\medskip
In fact, by (3.3), (3.9) \& $(\textrm{D9}_1)$, we have
\begin{equation*}
\begin{split}
\bigl[\,x_i^-&(0),\,x^-_{\alpha_{1,i+1}}(1)\,\bigr]_{r}\qquad (\hbox{by definition})\\
&=\bigl[\,x_i^-(0),\,\underbrace{\bigl[\,x_i^-(0),\,\bigl[\,x_{i-1}^-(0),\,x^-_{\alpha_{1,i-1}}(1)\,\bigr]_s\,\bigr]_s}\,\bigr]_r
\qquad(\hbox{using (3.3)})\\
&=\underbrace{\bigl[\,x_i^-(0),\,\bigl[\,\bigl[\,x_i^-(0),x_{i-1}^-(0)\,\bigr]_s,
\,x^-_{\alpha_{1,i-1}}(1)\,\bigr]_s\,\bigr]_r}\qquad(\hbox{this term using (3.3)})\\
&\quad+s\,\bigl[\,x_i^-(0),\,\bigl[\,x_{i-1}^-(0),\,\underbrace{\bigl[\,x_i^-(0),\,
x^-_{\alpha_{1,i-1}}(1)\,\bigr]}\,\bigr]\,\bigr]_r\qquad(=0\hbox{ by (3.5), (\textrm{D}$9_1$)})\\
&=\bigl[\,\underbrace{\bigl[\,x_i^-(0), x_i^-(0),
x_{i-1}^-(0)\,\bigr]_{(s,r)}},
\,x^-_{\alpha_{1,i-1}}(1)\,\bigr]_s\qquad(=0\hbox{ by (3.9)})\\
&\quad+\,r\,\bigl[\,\bigl[\,x_i^-(0),x_{i-1}^-(0)\,\bigr]_s,
\,\underbrace{\bigl[\,x_i^-(0),\,x^-_{\alpha_{1,i-1}}(1)\,\bigr]}\,\bigr]_{r^{-1}s}\qquad(=0\hbox{ by (3.5), (\textrm{D}$9_1$)}) \\
&=0.
\end{split}
\end{equation*}

\noindent {\bf Claim $(B)$}: \
$\bigl[\,x_i^-(0),\,x^-_{\alpha_{1,i+2}}(1)\,\bigr]_{\la\om_i',\om_{\alpha_{1,i{+}2}}\ra}=\bigl[\,x_i^-(0),\,x^-_{\alpha_{1,i+2}}(1)\,\bigr]=0$
\ ($i\ge 2$), \textit{if} $r\ne -s$.

\medskip
By definition, we note that
$\bigl[\,b,\,a\,\bigr]_u=-u\,\bigl[\,a,\,b\,\bigr]_{u^{-1}}$. So, we
get \begin{equation*}
\begin{split}\bigl[\,x_i^-(0),\,\underbrace{
x_{i+1}^-(0),\,x_i^-(0)\,}\bigr]_{(s,
r^{-1})}&=-s\,\bigl[\,x_i^-(0),\,
x_{i}^-(0),\,x_{i+1}^-(0)\,\bigr]_{(s^{-1}, r^{-1})}\quad(\hbox{by
 (3.6)})\\
&=-s\,\bigl[\,x_i^-(0),\,
x_{i}^-(0),\,x_{i+1}^-(0)\,\bigr]_{(r^{-1}, s^{-1})}\quad(\hbox{by
 (3.9)})\\
&=0.
\end{split}
\end{equation*}

We then consider the following deduction

\begin{equation*}
\begin{split}
\bigl[\,&x_i^-(0),\,x^-_{\alpha_{1,i+2}}(1)\,\bigr]_{r^{-1}s}=
\bigl[\,x_i^-(0),\,\underbrace{\bigl[\,x_{i+1}^-(0),\,\bigl[\,x_i^-(0),\,x^-_{\alpha_{1i}}(1)\,\bigr]_s\,\bigr]_s}\,\bigr]_{r^{-1}s}\
(\hbox{by
 (3.3)})\\
&=\underbrace{\bigl[\,x_i^-(0),\,\bigl[\,\bigl[\,x_{i+1}^-(0),\,x_i^-(0)\,\bigr]_s,\,x^-_{\alpha_{1i}}(1)\,\bigr]_s\,\bigr]_{r^{-1}s}}\qquad(\hbox{using
 (3.3)})\\
&\quad+s\bigl[\,x_i^-(0),\bigl[\,x_i^-(0),\underbrace{\bigl[\,x_{i+1}^-(0),\,x^-_{\alpha_{1i}}(1)\,\bigr]}\,\bigr]\,\bigr]_{r^{-1}s}\quad(\hbox{this
term$=0$ by (3.5),
(D9$_1$)})\\
&=\bigl[\,\underbrace{\bigl[\,x_i^-(0),\,x_{i+1}^-(0),\,x_i^-(0)\,\bigr]
_{(s,\,r^{-1})}},\,x^-_{\alpha_{1i}}(1)\,\bigr]_{s^2}\qquad(\hbox{this
term$=0$ by the above})\\
&\quad+\,r^{-1}\,\bigl[\,\bigl[\,x_{i+1}^-(0),\,x_i^-(0)\,\bigr]_s\,,\,
x^-_{\alpha_{1,i+1}}(1)\,\bigr]_{rs}\qquad(\hbox{using (3.4)})
\\
&=r^{-1}\,\bigl[\,\bigl[\,x_{i+1}^-(0),\,\underbrace{\bigl[\,x_i^-(0),\,
x^-_{\alpha_{1,i+1}}(1)\,\bigr]_r}\,\bigr]_{s^2}+\,\bigl[\,x^-_{\alpha_{1,i+2}}(1),\,
x_i^-(0)\,\bigr]_{r^{-1}s}\\
&=\bigl[\,x^-_{\alpha_{1,i+2}}(1),
x_i^-(0)\,\bigr]_{r^{-1}s}.\qquad(\hbox{$1$st term$=0$ by {\bf Claim
$(A)$}})
\end{split}
\end{equation*}
Expanding both sides of the above equation according to definition,
we easily get
$$(1+r^{-1}s)\,\bigl[\,x_i^-(0),\,x^-_{\alpha_{1,i+2}}(1)\,\bigr]=0.$$
Thus the required result is obtained under the assumption.

Now applying (3.5), we can get
\begin{equation*}
\begin{split}
[\,x_i^-(0),\,
x^-_{\theta}(1)\,]&=\bigl[\,x_{n-1}^-(0),\,\cdots,\,x_{i+2}^-(0),\,\underbrace{\bigl[\,x_i^-(0),\,
x^-_{\alpha_{1,i+2}}(1)\,\bigr]}\,\bigr]_{(s,\cdots,s)}\\
&=0.\qquad(\hbox{by {\bf Claim $(B)$}})
\end{split}
\end{equation*}

This completes the proof of Lemma 4.4.
\end{proof}

Next, we turn to check the relation below, whose argument (using the
quantum calculations) is crucial to our verification on
compatibilities of the defining relations system of $\mathcal
U_{r,s}(\widehat{\frak{sl}_n})$ mentioned in Remark 3.3.
\begin{prop} \
$[\,E_0, F_0\,]=\frac{\omega_0-\omega'_0}{r-s}$.
\end{prop}
\begin{proof} \
Using $(\textrm{D1})$ \& $(\textrm{D5})$, we have
\begin{equation*}
\begin{split}
\bigl[\,E_0,
F_0\,\bigr]&=\bigl[\,x^-_{\alpha_{1n}}(1)\,\gamma'^{-1}{\om_\theta}^{-1},\,\gamma^{-1}{\om'_\theta}^{-1}x^+_{\alpha_{1n}}(-1)\,\bigr]\\
&=\bigl[\,x^-_{\alpha_{1n}}(1),\,x^+_{\alpha_{1n}}(-1)\,\bigr]\,\cdot(\gamma^{-1}\gamma'^{-1}{\om_\theta}^{-1}{\om'_\theta}^{-1}).
\end{split}\tag{4.1}
\end{equation*}

Note that for $j\ge 1$, we have
\begin{gather*}
\bigl[\,x_{j+1}^-(0),
\,\om_j'\,\bigr]_s=(r-s)\,\om_j'\,x_{j+1}^-(0),\qquad
\bigl[\,x_{j+1}^-(0), \,\om_j\,\bigr]_s=0,\tag{4.2}\\
 \bigl[\,\om_j'\,X, x_{j+1}^+(k)\,\bigr]_r
=\om_j'\,\bigl[\,X, x_{j+1}^+(k)\,\bigr],\tag{4.3}\\
\bigl[\,x_j^+(k),
Y\om_{j+1}\,\bigr]_r=\bigl[\,x_j^+(k),Y\,\bigr]\,\om_{j+1}.\tag{4.4}
\end{gather*}
So (4.2) implies that there hold
\begin{gather*}
\bigl[\,x_{j+1}^-(0),
\,\bigl[\,x_j^-(0),\,x_j^+(0)\,\bigr]\,\bigr]_s=\om_j'\,x_{j+1}^-(0),\tag{4.5}\\
\bigl[\,x_2^-(0),
\,\bigl[\,x_1^-(1),\,x_1^+(-1)\,\bigr]\,\bigr]_s=\gamma\,\om_1'\,x_2^-(0),\tag{4.6}\\
\bigl[\,\bigl[\,x_{j+1}^-(0),\,x_{j+1}^+(0)\,\bigr],\,x_j^-(0)\,\big]_s=-x_j^-(0)\,\om_{j+1}.\tag{4.7}
\end{gather*}

Now let us write briefly
$$\bigl[\,x_1^+(-1),x_2^+(0),\cdots,x_{i-1}^+(0)\,\bigr]_{\la
r,\cdots,r\ra}:=\bigl[\,\bigl[\cdots\bigl[\,x_1^+(-1),x_2^+(0)\,\bigr]_r,\cdots\bigr]_r,x_{i-1}^+(0)\,\bigr]_r.$$
Thus, by (3.5), we have
\begin{equation*}
\begin{split}
\bigl[\,x^-_{\alpha_{1i}}(1),x^+_{\alpha_{1i}}(-1)\,\bigr]
&=\bigl[\,x^-_{\alpha_{1i}}(1),\,\bigl[\,x_1^+(-1),\,x_2^+(0),\cdots,x_{i-1}^+(0)\,\bigr]_{\la
r,\cdots,r\ra}\,\bigr]\\
&=\bigl[\,\underbrace{\bigl[\,x^-_{\alpha_{1i}}(1),
x_1^+(-1)\,\bigr]},\,x_2^+(0),\cdots,x_{i-1}^+(0)\,\bigr]_{\la
r,\cdots,r\ra}\\
&\quad
+\sum_{j=2}^{i-1}\,\bigl[\,x_1^+(-1),x_2^+(0),\cdots,\underbrace{\bigl[\,x^-_{\alpha_{1i}}(1),x_j^+(0)\,\bigr]},\cdots\,\bigr]_{\la
r,\cdots,r\ra}.\\
\end{split}\tag{4.8}
\end{equation*}

(i) \ For $j=1$, by (3.5),  (D8) \& (4.6), we have
\begin{equation*}
\begin{split}
\bigl[\,x^-_{\alpha_{1i}}(1),x_1^+(-1)\,\bigr]
&=\bigl[\,x_{i-1}^-(0), \cdots,\underbrace{\bigl[\,x_2^-(0),\bigl[\,x_1^-(1),x_1^+(-1)\,\bigr]\,\bigr]_s}\bigr]_{(s,\cdots,s)}\quad(\textrm{by (4.6)})\\
&=\gamma\om'_1\,x^-_{\alpha_{2i}}(0), \qquad (i>2)
\end{split}
\end{equation*}
so that
\begin{equation*}
\begin{split}
M(i):&=\bigl[\,\underbrace{\bigl[\,x^-_{\alpha_{1i}}(1),\,x_1^+(-1)\,\bigr]},x_2^+(0),\cdots\,x_{i-1}^+(0)\,\bigr]_{\la r,\cdots,r\ra}\\
&=\gamma\,\bigl[\,\underbrace{\bigl[\,\om_1'\,x^-_{\alpha_{2i}}(0),x_2^+(0)\,\bigr]_r},\cdots,x_{i-1}^+(0)\,\bigr]_{\la r,\cdots,r\ra}\\
&=\gamma\,\om_1'\,\bigl[\,\underbrace{\bigl[\,x^-_{\alpha_{2i}}(0),x_2^+(0)\,\bigr]},\cdots,x_{i-1}^+(0)\,\bigr]_{\la r,\cdots,r\ra}\\
&=\gamma\om_1'\,\bigl[\,\underbrace{\bigl[\,\om_2'\,x^-_{\alpha_{3i}}(0),x_3^+(0)\,\bigr]_r},\cdots,x_{i-1}^+(0)\,\bigr]_{\la r,\cdots,r\ra}\\
&=\gamma\om_1'\om_2'\,\bigl[\,\underbrace{\bigl[\,x^-_{\alpha_{3i}}(0),x_3^+(0)\,\bigr]},\cdots,x_{i-1}^+(0)\,\bigr]_{\la r,\cdots,r\ra}\\
&=\cdots\\
&=\gamma\om_1'\cdots\om_{i-2}'\,\bigl[\,x_{i-1}^-(0),\,x_{i-1}^+(0)\,\bigr]\\
&=\gamma\om_1'\cdots\om_{i-2}'\,\frac{\om_{i-1}'-\om_{i-1}}{r-s},
\qquad (i>2)
\end{split}
\end{equation*}
where we used the following identities, respectively
\begin{gather*}
\bigl[\,\om_{j-1}'\,x^-_{\alpha_{ji}}(0),x_j^+(0)\,\bigr]_r=\om_{j-1}'\,\bigl[\,x^-_{\alpha_{ji}}(0),x_j^+(0)\,\bigr],\qquad(\textrm{by (4.3)})\\
\bigl[\,x^-_{\alpha_{ji}}(0),x_j^+(0)\,\bigr]=\om_j'\,x^-_{\alpha_{j+1,i}}(0),\qquad(\textrm{by
(3.13) \& (4.5)}).
\end{gather*}

(ii) \ For $j=i-1$, again by (3.5), (3.3) \& (4.7), we get
\begin{equation*}
\begin{split}
\bigl[\,x^-_{\alpha_{1i}}(1),\,x_{i-1}^+(0)\,\bigr]&=\underbrace{\bigl[\,\bigl[\,x_{i-1}^-(0),\,x_{i-1}^+(0)\,\bigr],\,
\bigl[\,x_{i-2}^-(0),\,x_{\alpha_{1,i-2}}^-(1)\,\bigr]_s\,\bigr]_s}\quad(\textrm{by (3.3)})\\
&=\bigl[\,\underbrace{\bigl[\,\bigl[\,x_{i-1}^-(0),\,x_{i-1}^+(0)\,\bigr],\,x_{i-2}^-(0)\,\bigr]_s},\,x^-_{\alpha_{1,i-2}}(1)\,\bigr]_s
\quad(\textrm{by (4.7)})\\
&\quad+s\,\bigl[\,x_{i-2}^-(0),\,\underbrace{\bigl[\,\bigl[\,x_{i-1}^-(0),\,x_{i-1}^+(0)\,\bigr],\,x^-_{\alpha_{1,i-2}}(1)\,\bigr]}\,\bigr]\quad(\textrm{$=0$})\\
&=-\bigl[\,x_{i-2}^-(0)\,\om_{i-1},\,x^-_{\alpha_{1,i{-}2}}(1)\,\bigr]_s\\
&=-\bigl[\,x_{i-2}^-(0),\,x^-_{\alpha_{1,i{-}2}}(1)\,\bigr]_s\,\om_{i-1}\\
&=-\,x^-_{\alpha_{1,i{-}1}}(1)\,\om_{i-1},
\end{split}
\end{equation*}
where we notice that (*):
$\bigl[\,\bigl[\,x_{i-1}^-(0),\,x_{i-1}^+(0)\,\bigr],\,x^-_{\alpha_{1,i-2}}(1)\,\bigr]=0$.

\medskip
Thereby, we further obtain
\begin{equation*}
\begin{split}
N(i):&=\bigl[\,x_1^+(-1),\,x_2^+(0),\cdots,\underbrace{\bigl[\,x^-_{\alpha_{1i}}(1),x_{i-1}^+(0)\,\bigr]},\cdots\,\bigr]_{\la
r,\cdots,r\ra}\\
&=-\underbrace{\bigl[\,x_{\alpha_{1,i-1}}^+(-1),
x^-_{\alpha_{1,i{-}1}}(1)\,\om_{i-1}\,\bigr]_r}\qquad(\textrm{by (4.4)})\\
&=-\bigl[\,x_{\alpha_{1,i-1}}^+(-1),
x^-_{\alpha_{1,i{-}1}}(1)\,\bigr]\,\om_{i-1}\\
&=\bigl[\,x^-_{\alpha_{1,i{-}1}}(1),
x_{\alpha_{1,i-1}}^+(-1)\,\bigr]\,\om_{i-1}.
\end{split}
\end{equation*}

\medskip
(iii) \ For $1<j<i-1$, by (3.5), (3.3), (4.7) \& (D9$_1$), we
obtain
\begin{equation*}
\begin{split}
\bigl[\,x^-_{\alpha_{1i}}(1),x_j^+(0)\,\bigr]&=\bigl[\,x_{i-1}^-(0),\cdots
\underbrace{\bigl[\,\bigl[\,x_j^-(0),x_j^+(0)\,\bigr],\bigl[\,x_{j-1}^-(0),x_{\alpha_{1,j-1}}^-(1)\,\bigr]_s\bigr]_s}\,\bigr]_{(s,\cdots,s)}\\
&\hskip4.8cm\quad(\textrm{by (3.3)})\\
&=\bigl[\,x_{i{-}1}^-(0),\cdots
\underbrace{\bigl[\,\bigl[\,x_j^-(0),x_j^+(0)\,\bigr],x_{j-1}^-(0)\,\bigr]_s},x_{\alpha_{1,j-1}}^-(1)\,\bigr]_{(s,\cdots,s)}\\
&\hskip3.5cm\quad(\textrm{by (4.7)})\\
&\quad+s\,\bigl[\,x_{i{-}1}^-(0),\cdots\bigl[\,x_{j-1}^-(0),
\underbrace{\bigl[\bigl[x_j^-(0),x_j^+(0)\bigr],x_{\alpha_{1,j-1}}^-(1)\,\bigr]}\,\bigr]\bigr]_{(s,\cdots,s)}\\
&\hskip6cm(\textrm{$=0$ by (*)})\\
&=-\,\bigl[\,x_i^-(0),\cdots,x_{j+1}^-(0),\,\underbrace{\bigl[\,
x_{j-1}^-(0)\,\om_j,\,x_{\alpha_{1,j-1}}^-(1)\,\bigr]_s}\,\bigr]_{(s,\cdots,s)}\\
&=-\,\bigl[\,x_i^-(0),\,\cdots,\underbrace{\bigl[\,x_{j+1}^-(0),\,x_{\alpha_{1j}}^-(1)\,\om_j\,\bigr]_s
}\,\bigr]_{(s,\cdots,s)}\\
&=-\,\bigl[\,x_i^-(0),\cdots\underbrace{\bigl[x_{j+1}^-(0),x_{\alpha_{1j}}^-(1)\bigr]}\om_j\,\bigr]_{(s,
\cdots,s)}\quad(\textrm{by (3.5), (D9$_1$)})\\
&=0,
\end{split}
\end{equation*}
where in forth and fifth equality ``=" we used the following
identities, respectively
\begin{gather*}
\bigl[\,x_{j-1}^-(0)\,\om_j,\,x_{\alpha_{1,j-1}}^-(1)\,\bigr]_s=\bigl[\,x_{j-1}^-(0),\,x_{\alpha_{1,j-1}}^-(1)\,\bigr]_s\,\om_j
=x_{\alpha_{1j}}^-(1)\,\om_j,\\
\bigl[\,x_{j+1}^-(0),\,x_{\alpha_{1j}}^-(1)\,\om_j\,\bigr]_s=\bigl[\,x_{j+1}^-(0),\,x_{\alpha_{1j}}^-(1)\,\bigr]\,\om_j.
\end{gather*}

\medskip
As a result of (i), (ii) \& (iii),  (4.8) becomes
\begin{equation*}
\begin{split}
\bigl[\,x^-_{\alpha_{1i}}(1), x^+_{\alpha_{1i}}(-1)\,\bigr]
&=M(i)+N(i)\\
&=M(i)+\bigl[\,x^-_{\alpha_{1,i{-}1}}(1),
x_{\alpha_{1,i-1}}^+(-1)\,\bigr]\,\om_{i-1}\\
&=M(i)+M(i{-}1)\,\om_{i-1}+\bigl[\,x^-_{\alpha_{1,i{-}2}}(1),
x_{\alpha_{1,i-2}}^+(-1)\,\bigr]\,\om_{i-2}\om_{i-1}
\\
&=\cdots\\
&=M(i)+M(i{-}1)\,\om_{i-1}+M(i{-}2)\,\om_{i-2}\om_{i-1}+\cdots\\
&\quad+M(3)\,\om_3\cdots\om_{i-1}+\bigl[\,x^-_{\alpha_{12}}(1),
x_{\alpha_{12}}^+(-1)\,\bigr]\,\om_2\cdots\om_{i-1}\\
&=\frac{\gamma\om_{\alpha_{1i}}'-\gamma'\om_{\alpha_{1i}}}{r-s},
\qquad (i>1),
\end{split}\tag{4.9}
\end{equation*}
where we used (D8) to get
$$\bigl[\,x^-_{\alpha_{12}}(1),
x_{\alpha_{12}}^+(-1)\,\bigr]=\frac{\gamma\om_1'-\gamma'\om_1}{r-s}.$$

\medskip
Therefore, by (4.9), (4.1) takes the required formula:
$$
[\,E_0,\,F_0\,]=\frac{\gamma'^{-1}\,\om_\theta^{-1}-\gamma^{-1}\,\om_\theta'^{-1}}{r-s}.
$$

The proof of Proposition 4.5 is complete.
\end{proof}

For $(\textrm{A5})$: \ we need only to verify that $[\,E_0,
E_j\,]=0$ and $[\,F_0,\,F_j\,]=0$ for $1<j<n-1$. Actually, in the
proof of Proposition 4.5, the fact that
$\bigl[\,x^-_{\alpha_{1i}}(1), x_j^+(0)\,\bigr]=0$ for $1<j<i-1$
implies the first identity (taking $i=n$) since $[\,E_0,
E_j\,]=[\,x^-_{\theta}(1),x^+_j(0)\,]$
$\cdot\,\gamma'^{-1}\om_\theta^{-1}$. The second can be obtained
utilizing $\tau$ on the first one.\hfill$\Box$

\medskip
For $(\textrm{A6})$: \ when $i\cdot j \neq 0$, $(\textrm{D9$_n$})$
implies that the corresponding generators satisfy exactly those
$(r,s)$-Serre relations in $U_{r,s}(\frak {sl}_n)$. So, it is enough
to check the $(r,s)$-Serre relations involving the indices with
$i\cdot j=0$.

\begin{lemm} \ $(1)$ \ $E_0E_1^2-(r+s)E_1E_0E_1+(r
s)\,E_1^2E_0=0$,

$(2)$ \ $E_0^2E_1-(r+s)E_0E_1E_0+(rs)\,E_1E_0^2=0$,

$(3)$ \
$E_{n-1}^2E_0-(r+s)\,E_{n-1}E_0E_{n-1}+(rs)\,E_0E_{n-1}^2=0$,

$(4)$ \ $E_{n-1}E_0^2-(r+s)\,E_0E_{n-1}E_0+(rs)\,E_0^2E_{n-1}=0$,

$(5)$ \ $F_1^2F_0-(r+s)F_1F_0F_1+(r s)\,F_0F_1^2=0$,

$(6)$ \ $F_1F_0^2-(r+s)F_0F_1F_0+(rs)\,F_0^2F_1=0$,

$(7)$ \
$F_0F_{n-1}^2-(r+s)\,F_{n-1}F_0F_{n-1}+(rs)\,F_{n-1}^2F_0=0$,

$(8)$ \ $F_0^2F_{n-1}-(r+s)\,F_0F_{n-1}F_0+(rs)\,F_{n-1}F_0^2=0$.
\end{lemm}
\begin{proof} \ The proofs for the relations of $(5)$---$(8)$ follow from taking $\tau$ on
the first four relations $(1)$---$(4)$. We shall demonstrate the
first two $(r,s)$-Serre relations, the third and forth ones are
similar to the first two relations (1) \& (2), which are left to the
reader.

(1) \  Observing
\begin{equation*}
\begin{split}
\bigl[\,E_1, x^-_{\theta}(1)\,\bigr]&=\bigl[\,x_{n-1}^-(0),
\cdots,\,\underbrace{\bigl[\,x_1^+(0),\,
x_1^-(1)\,\bigr]}\,\bigr]_{(s,\cdots,s)}
\;\quad\qquad\qquad(\hbox{using $(\textrm{D8})$})\\
&=\gamma^{-\frac{1}2}\,\bigl[\,x_{n-1}^-(0),\cdots,\,\underbrace{\bigl[\,x_2^-(0),\,\om_1\,a_1(1)\,\bigr]_s}\,\bigr]_{(s,\cdots,s)}
\qquad(\hbox{using $(\textrm{D5})$})\\
&=\gamma^{-\frac{1}2}\,s\,\om_1\,\bigl[\,x_{n-1}^-(0),\cdots,\,\underbrace{\bigl[\,x_2^-(0),\,a_1(1)\,\bigr]}\,\bigr]_{(s,\cdots,s)}
\qquad(\hbox{using $(\textrm{D6})$})\\
&=-(rs)^{-\frac1{2}}x_{\alpha_{2n}}^-(1)\,\om_1,
\end{split}
\end{equation*}
we have
\begin{equation*}
\begin{split}
 E_0&E_1^2-(r+s)E_1E_0E_1+(r s)\,E_1^2E_0 \\
&=(r s)\Big(E_1^2 x^-_{\theta}(1)-(1+r^{-1}s)E_1
x^-_{\theta}(1)E_1+(r^{-1}s)\,x^-_{\theta}(1)E_1^2\Big)
(\gamma'^{-1}\om_{\theta}^{-1})\\
&= (rs)\,\bigl[\,E_1, \,\underbrace{\bigl[\,E_1,\,x^-_{\theta}(1)\,\bigr]}\,\bigr]_{r^{-1}s}\,(\gamma'^{-1}\om_{\theta}^{-1})\qquad(\hbox{by (3.6)})\\
&=-(rs)^{\frac1{2}}\,\underbrace{\bigl[\,x_1^+(0),\,
x_{\alpha_{2n}}^-(1)\,\om_1\,\bigr]_{r^{-1}s}}\, (\gamma'^{-1}\om_{\theta}^{-1})\\
&=-(rs)^{\frac1{2}}\,\,\underbrace{\bigl[\,x_1^+(0),
x_{\alpha_{2n}}^-(1)\,\bigr]}\,\om_1\, (\gamma'^{-1}\om_{\theta}^{-1})\\
&=0.\qquad(\hbox{by (3.5), $(\textrm{D8})$})
\end{split}
\end{equation*}

(2) \ Using the formula of $\,\bigl[\,E_1,
x^-_{\theta}(1)\,\bigr]\,$ derived in (1) above, we have
\begin{equation*}
\begin{split}
E_0^2&E_1-(r+s)E_0E_1E_0+(rs)\,E_1E_0^2 \\
&=(rs)\,\bigl[\,x^-_{\theta}(1), \underbrace{x^-_{\theta}(1),
E_1}\,\bigr]_{(1,\,rs^{-1})}(\gamma'^{-2} \om_{\theta}^{-2})\\
&=(rs)^{\frac{1}2}\,\bigl[\,x^-_{\theta}(1),\,
x_{\alpha_{2n}}^-(1)\,\om_1\,\bigr]_{rs^{-1}}\,(\gamma'^{-2} \om_{\theta}^{-2}) \\
&=(rs)^{\frac{1}2}\,\bigl[\,x^-_{\theta}(1),\,
x_{\alpha_{2n}}^-(1)\,\bigr]_{s^{-1}}\,\om_1\,(\gamma'^{-2} \om_{\theta}^{-2})\\
&=-(rs^{-1})^{\frac{1}2}\,\underbrace{\bigl[\,x_{\alpha_{2n}}^-(1),
x^-_{\theta}(1)\,\bigr]_s}\,\om_1\,(\gamma'^{-2}\,\om_{\theta}^{-2})\\
&=0,\qquad(\textrm{by {\bf Claim} ($C$) below})
\end{split}
\end{equation*}
where we used the following claim:

\medskip
\noindent {\bf Claim $(C)$}: \ $\bigl[\,x_{\alpha_{2n}}^-(1),
x^-_{\alpha_{1n}}(1)\,\bigr]_s=0$, \textit{for}  $n>2$ \textit{and}
$r\ne -s$.
\medskip

The argument for Claim $(C)$ is technical. Indeed, by induction on
$n$, we have:

When $n=3$, by (3.8), one gets
\begin{equation*}
\begin{split}
\bigl[\,x_2^-(1),\,x^-_{\alpha_{13}}(1)\,\bigr]_s&=\bigl[\,x_2^-(1),\,\underbrace{\bigl[\,x_2^-(0),\,x_1^-(1)\,\bigr]_s}\,\bigr]_s\qquad(\textrm{by (3.8)})\\
&=-(rs)^{\frac1{2}}\,\bigl[\,x_2^-(1),\,\underbrace{\bigl[\,x_1^-(0),\,x_2^-(1)\,\bigr]_{r^{-1}}}\,\bigr]_s\\
&=(rs^{-1})^{-\frac1{2}}\,\bigl[\,x_2^-(1),\,x_2^-(1),\,x_1^-(0)\,\bigr]_{(r,
s)}\\
&=(rs^{-1})^{-\frac1{2}}\,\bigl[\,x_2^-(1),\,x_2^-(1),\,x_1^-(0)\,\bigr]_{(s,
r)}\qquad(\textrm{by (3.6)})\\& =0,\qquad(\textrm{by (3.9)})
\end{split}
\end{equation*}
which is exactly the $(r,s)$-Serre relation (see (3.9)).

While for $n>3$, we first notice the fact:
$$
\bigl[\,x_{n-1}^-(0),\, x_{\alpha_{2n}}^-(1)\,\bigr]_{\la
\om_{n-1}',\om_{\alpha_{2n}}\ra}=\bigl[\,x_{n-1}^-(0),\,
x_{\alpha_{2n}}^-(1)\,\bigr]_r=0, \qquad\hbox{\it for } \
n>3,\leqno(4.10)
$$
which can be proved using the same method of the proof of (II) in
Lemma 4.4.

We thus have
\begin{equation*}
\begin{split}
\bigl[\,x_{\alpha_{2n}}^-(1),&\,
x^-_{\theta}(1)\,\bigr]_r=\bigl[\,\bigl[\,x_{n-1}^-(0),\,x_{\alpha_{2,n-1}}^-(1)\,\bigr]_s,\,
x^-_{\theta}(1)\,\bigr]_r\qquad(\hbox{by (3.4)})\\
&=\bigl[\,x_{n-1}^-(0),\,\bigl[\,x_{\alpha_{2,n-1}}^-(1),\,
x^-_{\theta}(1)\,\bigr]_1\,\bigr]_{r
s}\\
&\quad+\bigl[\,\underbrace{\bigl[\,x_{n-1}^-(0),\,
x^-_{\theta}(1)\,]_r},\, x_{\alpha_{2,n-1}}^-(1)
\,\bigr]_s\qquad(\hbox{$=0$ by {\bf Claim} (A)})
\\
&=\bigl[\,x_{n-1}^-(0),\,\underbrace{\bigl[\,x_{\alpha_{2,n-1}}^-(1),\,
\bigl[\,x_{n-1}^-(0),\,x^-_{\alpha_{1,n-1}}(1)\,\bigr]_{s}\,\bigr]_1}\,\bigr]_{r s}\qquad (\hbox{by (3.3)})\\
&=\bigl[\,x_{n-1}^-(0),\bigl[\,\underbrace{\bigl[\,x_{\alpha_{2,n-1}}^-(1),\,
x_{n-1}^-(0)\,\bigr]_{s^{-1}}},\,x_{\alpha_{1,n-1}}^-(1)
\,\bigr]_{s^2}\,\bigr]_{r s}\\
&\quad+s^{-1}\bigl[\,x_{n-1}^-(0),\,\bigl[\,x_{n-1}^-(0),\,\underbrace{\bigl[\,x_{\alpha_{2,n-1}}^-(1),\,
x_{\alpha_{1,n-1}}^-(1)\,\bigr]_{s}}\,\bigr]_{s^2}\,\bigr]_{r s}\\
&\hskip2cm(\hbox{$2$nd sumand$=0$ using induction
hypothesis})\\
&=-s^{-1}\bigl[\,x_{n-1}^-(0),\,\bigl[\,x_{\alpha_{2n}}^-(1),
\,x_{\alpha_{1,n-1}}^-(1)\,\bigr]_{s^2}\,\bigr]_{r s}\qquad(\hbox{by
(3.3)}) \\
\end{split}
\end{equation*}
\begin{equation*}
\begin{split}
&=-s^{-1}\bigl[\,\underbrace{\bigl[\,x_{n-1}^-(0),\,x_{\alpha_{2n}}^-(1)\,\bigr]_r},
\,x_{\alpha_{1,n-1}}^-(1)\,\bigr]_{s^3}\qquad( \hbox{$=0$
by (4.10)}) \\
&\quad-rs^{-1}\bigl[\,x_{\alpha_{2n}}^-(1),
\,\underbrace{\bigl[\,x_{n-1}^-(0),\,x_{\alpha_{1,n-1}}^-(1)\,\bigr]_s}\,\bigr]_{r^{-1}s^2}\qquad (\hbox{by definition})\\
&=-rs^{-1}\bigl[\,x_{\alpha_{2n}}^-(1),
\,x^-_{\theta}(1)\,\bigr]_{r^{-1}s^2}.
\end{split}
\end{equation*}
By definition, expanding both sides of the above identity gives us
$$(1+rs^{-1})\,x_{\alpha_{2n}}^-(1)\,x_{\alpha_{1n}}^-(1)=(r+s)\,x_{\alpha_{1n}}^-(1)\,x_{\alpha_{2n}}^-(1),$$
which means $\bigl[\,x_{\alpha_{2n}}^-(1),\,
x_{\alpha_{1n}}^-(1)\,\bigr]_s=0$, under the assumption $r\ne -s$.
\medskip

For $(\textrm{A7})$: \ the verification is analogous to that of
$(\textrm{A6})$.
\end{proof}

\noindent{\bf 4.3.} \ {\it Proof of Theorem 4.2.} \ We shall show
that the algebra ${\mathcal U}_{r,s}(\widehat{\frak {sl}_n})$ is
generated by $E_i,\, F_i,\,\om_i^{\pm1},\,{\om}_i'^{\,\pm1}$,
$\gamma^{\pm\frac{1}2}$, $\gamma'^{\,\pm\frac{1}2},\, D^{\pm1},\,
D'^{\,\pm1}$ ($i\in I_0$).

To this end, we need to prove the following results.
\begin{lemm} \ $(1)$ \
$x_1^-(1)=\bigl[\,E_2,\,E_3,\cdots,\,E_{n-1},
E_0\,\bigr]_{(r,\cdots,r)} \gamma'\om_1\in {\mathcal
U\,'}_{r,s}(\widehat{\frak {sl}_n})$, then for any $i\in I$,
$x_i^-(1)\in {\mathcal U\,'}_{r,s}(\widehat{\frak {sl}_n})$.

$(2)$ \ $x_1^+(-1)=\tau\Bigl(\bigl[E_2,E_3,\cdots,E_{n-1},
E_0\bigr]_{(r,\cdots,r)} \gamma'\om_1\Bigr)=\gamma\om'_1\bigl[\,F_0,
F_{n-1}, \cdots,F_3$, $F_2\,\bigr]_{\la s,\cdots, s\ra}\in {\mathcal
U\,'}_{r,s}(\widehat{\frak {sl}_n})$, then for any $i\in I$,
$x_i^+(-1)\in {\mathcal U\,'}_{r,s}(\widehat{\frak {sl}_n})$.
\end{lemm}
\begin{proof} \ $(1)$ \ Set
$\widetilde{E}(i)=x_{\alpha_{1,i{+}1}}^-(1)\,\om_{i+1}\cdots\om_{n-1}\gamma'^{-1}
\om_{\theta}^{-1}$ for $i\ge1$, where $\widetilde{E}(n{-}1)=E_0$.
Observing that
$\bigl[\,x_i^+(0),\,x_{\alpha_{1,i{+}1}}^-(1)\,\bigr]=x_{\alpha_{1i}}^-(1)\,\om_i$
in the proof (see, case (ii)) of Proposition 4.5, we get an
important recursive relation:
\begin{equation*}
\begin{split}
 \bigl[\,E_i, \, \widetilde{E}(i)\,\bigr]_r
 &=\bigl[\,x^+_i(0),\,x_{\alpha_{1,i{+}1}}^-(1) \,\om_{i+1}\cdots\om_{n-1}\gamma'^{-1} \om_{\theta}^{-1}\, \bigr]_r\\
&=\bigl[\,x_i^+(0),\,x_{\alpha_{1,i{+}1}}^-(1)\,\bigr]\,\om_{i+1}\cdots\om_{n-1}\gamma'^{-1}\om_{\theta}^{-1}\\
&=\widetilde{E}(i{-}1).
\end{split}\tag{4.11}
\end{equation*}
Recursively using the above relations, we obtain
\begin{equation*}
\begin{split}
x_1^-(1)&=\widetilde{E}(1)\gamma'\om_1=\bigl[\,E_2, \, \widetilde{E}(2)\,\bigr]_r\gamma'\om_1=\cdots\\
&=\bigl[\,E_2,\cdots,E_{n-1},
\widetilde{E}(n{-}1)\,\bigl]_{(r,\cdots,r)}\,\gamma'\om_1\\
&=\bigl[\,E_2,\cdots,E_{n-1},\,E_0\,\bigl]_{(r,\cdots,r)}\,\gamma'\om_1\\
&\in {\mathcal U\,'}_{r,s}(\widehat{\frak {sl}_n}).
\end{split}\tag{4.12}
\end{equation*}

Now suppose that we already have obtained $x_i^-(1)\in{\mathcal
U\,'}_{r,s}(\widehat{\frak {sl}_n})$ for $i\ge 1$. Notice that
\begin{equation*}
\begin{split}
x_{i+1}^-(1)&=(rs)\,\underbrace{\bigl[\,\bigl[\,x_i^+(0),\,x_i^-(0)\,\bigr],\,x_{i+1}^-(1)\,\bigr]_{r^{-1}}}\,\om_i^{-1}\qquad\hbox{(by
(3.4))}\\
&=(rs)\,\bigl[\,x_i^+(0),\,\underbrace{x_i^-(0),\,x_{i+1}^-(1)}\,\bigr]_{(r^{-1},1)}\,\om_i^{-1}\qquad\;\hbox{(by
(3.8))}\\
&=-(rs)^{\frac1{2}}\bigl[\,x_i^+(0),\,x_{i+1}^-(0),\,x_i^-(1)\,\bigr]_{(s,1)}\,\om_i^{-1}\\
&=(rs)^{\frac1{2}}\bigl[\,\bigl[\,F_{i+1},\,x_i^-(1)\,\bigr]_s,\,E_i\,\bigr]\,\om_i^{-1}\\
&\in {\mathcal U\,'}_{r,s}(\widehat{\frak {sl}_n}),
\end{split}\tag{4.13}
\end{equation*}
which gives rise to the recursive construction of some basic quantum
real root vectors of level $1$. Hence, we obtain the required
result.

\medskip
$(2)$ \ Set 
$\widetilde{F}(i)=\tau(\widetilde{E}(i))=\gamma^{-1}\om_{\theta}'^{-1}\om_{n-1}'\cdots\om_{i+1}'\,
x_{\alpha_{1,i{+}1}}^+(-1)$ for $i\ge1$, where
$\widetilde{F}(n{-}1)=F_0$. Applying $\tau$ to (4.11), we see that
$\bigl[\,\widetilde{F}(i),\,F_i\,\bigr]_s=\widetilde{F}(i{-}1)$
and
$\widetilde{F}(1)=
\gamma^{-1}\om_1'^{-1}x_1^+(-1)$, which implies the first claim.

The remaining claim follows from
$$
x_{i+1}^+(-1)=\tau(x_{i+1}^-(1))=(rs)^{\frac1{2}}\om_i'^{-1}\,\bigl[\,F_i,\,\bigl[\,x_i^+(-1),\,E_{i+1}\,\bigr]_r\,\bigr]
\in {\mathcal U\,'}_{r,s}(\widehat{\frak {sl}_n}).\leqno(4.14)$$

This completes the proof of Lemma 4.7.
\end{proof}

We observe that Lemma 4.7, together with (4.12), (4.13) \& (4.14),
gives the construction of the Drinfel'd generators of level $1$.
Furthermore, the first conclusion of the following Lemma gives the
the construction of the quantum imaginary root vectors of any level
($\ne0$), while the second gives the construction of some basic
quantum real root vectors of any level.

Actually, as a result of Definition 3.9 and Lemma 4.8 below, this
approach also gives the construction of all quantum real root
vectors of any level.
\begin{lemm} \ $(1)$ \ $a_i(\ell)\in {\mathcal U\,'}_{r,s}(\widehat{\frak {sl}_n})$, \ for $\;\ell\in\mathbb{Z}\backslash \{0\}$.

$(2)$ \ $x_i^{\pm}(k)\in {\mathcal U\,'}_{r,s}(\widehat{\frak
{sl}_n})$, \ for $\;k\in\mathbb{Z}$.
\end{lemm}
\begin{proof} \ $(1)$ \ At first, it follows from $(\textrm{D8})$
that
\begin{gather*}
a_i(1)=\om_i^{-1}\gamma^{1/2}\,\bigl[\,x_i^+(0),\,x_i^-(1)\,\bigr]\in {\mathcal U\,'}_{r,s}(\widehat{\frak {sl}_n}),\tag{4.15}\\
a_i(-1)=\om_i'^{-1}{\gamma}'^{1/2}\,\bigl[\,x_i^+(-1),\,x_i^-(0)\,\bigr]=\tau\bigl(a_i(1)\bigr)\in
{\mathcal U\,'}_{r,s}(\widehat{\frak {sl}_n}).\tag{4.16}
\end{gather*}

Suppose that we have already obtained $a_i(\pm\, \ell')\in {\mathcal
U\,'}_{r,s}(\widehat{\frak {sl}_n})\;$ for all $\;\ell'\le \ell$ and
some $\;\ell\ge 1$.

Now using $(\textrm{D6$_n$})$ \& $(\textrm{D8})$, we have the
following expansion (in fact, the expansions of both sides are the
same which also show the compatibility between (D6$_n$) and (D8) for
$n=1, 2$.):
\begin{equation*}
\begin{split}
{\mathcal U\,'}_{r,s}(\widehat{\frak {sl}_n})
&\ni\bigl[\,x_i^+(0),\,\bigl[\,a_i(\ell),\, x_i^-(1)\,\bigr]\,\bigr]\\
&=\bigl[\,\bigl[\,x_i^+(0),\, a_i(\ell)\,\bigr], \,x_i^-(1)\,\bigr]+\bigl[\,a_i(\ell),\, \bigl[\,x_i^+(0),\, x_i^-(1)\,\bigr]\,\bigr]\\
&=\,*\gamma'^{\frac{\ell}{2}}\,\bigl[\,x_i^+(\ell),\,
x_i^-(1)\,\bigr]\\
&\quad+\bigl[\,a_i(\ell),\,
\gamma^{-\frac{1}{2}}\,\om_i\,a_i(1)\,\bigr]\quad(\textrm{this term$=0$ by (D2)})\\
&=*(\gamma\gamma')^{-\frac{\ell}{2}}\gamma^{-\frac{1}2}\,\om_i\,\left[a_i(\ell{+}1)+\sum_{1<p\le
\ell{+}1\atop\sum_k
\ell_k=\ell{+}1}*'(r{-}s)^{p-1}a_i(\ell_{j_1})\cdots
a_i(\ell_{j_p})\right],
\end{split}\tag{4.17}
\end{equation*}
where scalars $*,\,*'\in\mathbb{K}\backslash \{0\}$. So
$a_i(\ell{+}1)\in {\mathcal U\,'}_{r,s}(\widehat{\frak {sl}_n})$.

Applying $\tau$ to the above formula, we can get
$a_i(-(\ell{+}1))\in {\mathcal U\,'}_{r,s}(\widehat{\frak {sl}_n})$.
Thereby, $a_i(\ell)\in {\mathcal U\,'}_{r,s}(\widehat{\frak
{sl}_n})$, for any $\ell\in\mathbb{Z}\backslash\{0\}$.

$(2)$ follows from $(\textrm{D6})$ (setting $i=j$ and $k=0$),
together with $(1)$.
\end{proof}

Therefore, we have proved ${\mathcal U\,'}_{r,s}(\widehat{\frak
{sl}_n})={\mathcal U}_{r,s}(\widehat{\frak {sl}_n})$, that is to
say, the latter is indeed generated by $E_i,\,
F_i,\,\om_i^{\pm1},\,{\om}_i'^{\,\pm1}$, $\gamma^{\pm\frac{1}2}$,
$\gamma'^{\,\pm\frac{1}2},\, D^{\pm1},\, D'^{\,\pm1}$ ($i\in
I_0$).\hfill$\Box$

\medskip
\noindent{\bf 4.4} \ {\it Proof of Theorem 4.3.} \  From {\bf 4.2}
\& {\bf 4.3}, we actually get an algebra epimorphism $\Psi:\;
U_{r,s}(\widehat{\frak {sl}_n}) \longrightarrow {\cal
U}_{r,s}(\widehat{\frak {sl}_n})$, since both algebras have the
essentially same generators system enjoying with the defining
relations from the former.

Notice that both algebras $U_{r,s}(\widehat{\frak {sl}_n})$ and
${\cal U}_{r,s}(\widehat{\frak {sl}_n})$ have commonly a natural
$Q$-gradation structure (see Corollary 2.8), which is by definition
preserved evidently under $\Psi$. On the other hand, both toral
subalgebras $U_{r,s}(\widehat{\frak {sl}_n})^0$ and ${\cal
U}_{r,s}(\widehat{\frak {sl}_n})^0$ generated by the same generators
system of group-like elements
$$\{\,\om_i^{\pm1},
\,{\om'}_i^{\pm1}\, (i\in I_0),\,\gamma^{\pm\frac{1}2},
\gamma'^{\,\pm\frac{1}2},\,D^{\pm1},\, D'^{\,\pm1}\,\}$$ are
obviously isomorphic with respect to
$\Psi^0:=\Psi|_{U_{r,s}(\widehat{\frak {sl}_n})^0}$.

Assigned to the positive or negative nilpotent Lie subalgebra
$\widehat {\frak n}^\pm$ of $\widehat{\frak {sl}_n}$ are two
subalgebras $U_{r,s}(\widehat{\frak n}^\pm)$ and ${\cal
U}_{r,s}(\widehat{\frak n}^\pm)$. Both are generated by
$\widehat{\frak n}^\pm$ in $U_{r,s}(\widehat{\frak {sl}_n})$ and
${\cal U}_{r,s}(\widehat{\frak {sl}_n})$ respectively. Denote
$\Psi^\pm:=\Psi|_{U_{r,s}(\widehat{\frak n}^\pm)}$. By Corollary
2.7, the double structure of $U_{r,s}(\widehat{\frak {sl}_n})$ in
Theorem 2.5 implies its triangular decomposition structure
$U_{r,s}(\widehat{\frak n}^-)\otimes U_{r,s}(\widehat{\frak
{sl}_n})^0\otimes U_{r,s}(\widehat{\frak n}^+)$. This fact likewise
indicates that $\Psi$ has a corresponding decomposition
$\Psi^-\otimes\Psi^0\otimes \Psi^+$. So, we are left to show
$\Psi^\pm$ are isomorphic. It suffices to consider the epimorphism
$\Psi^+:\;U_{r,s}(\widehat{\frak n}^+)\longrightarrow {\cal
U}_{r,s}(\widehat{\frak n}^+)$.

Observe that $U_{r,s}(\widehat{\frak n}^+)$ (resp. ${\cal
U}_{r,s}(\widehat{\frak n}^+)$) is generated by elements $e_i$
(resp. $E_i$) for $i\in I_0$ and subject to $(r,s)$-Serre relations
$(A5)$ \& $(A6)$. To check that $\Psi^+$ is an isomorphism, now we
fix $r=q$ and specialize $s$ at $q^{-1}$ as follows.

Note that $U_{r,s}(\widehat{\frak{n}}^+)$ can be viewed as to be
defined over the Laurent polynomials ring
$\mathbb{Q}[r^{\pm1},s^{\pm1}]$. Let ${\cal A}\subset
\mathbb{Q}(r,s)$ be the localization of ring
$\mathbb{Q}[r^{\pm1},s^{\pm1}]$ at the maximal ideal $(rs-1)$. Let
$U_{\cal A}^+$ be the ${\cal A}$-subalgebra of
$U_{r,s}(\widehat{\frak{n}}^+)$ generated by $e_i$ $(i\in I_0)$. Let
$(rs-1)U_{\cal A}^+$ be the ideal generated by $(rs-1)$ in $U_{\cal
A}^+$. Define the algebra $U_q^+$, the specialization of
$U_{r,s}(\widehat{\frak{n}}^+)$ at $s=q^{-1}$, by $U_q^+=U_{\cal
A}^+/(rs-1)U_{\cal A}^+$. Obviously, $U_q^+\cong
U_q(\widehat{\frak{n}}^+)$, the usual one-parameter quantum
subalgebra of $U_q(\widehat{\frak {sl}_n})$. However, in this case,
$\Psi^+$ induces the isomorphism $\Psi^+: \;U_q(\widehat{\frak
n}^+)\longrightarrow {\cal U}_q(\widehat{\frak n}^+)$ given by the
Drinfel'd isomorphism in the one-parameter case (see [B2] or [J2]).

Since specialization doesn't change the root multiplicities,
$\Psi^+:\; U_{r,s}(\widehat{\frak n}^+)\longrightarrow {\cal
U}_{r,s}(\widehat{\frak n}^+)$ is an isomorphism.\hfill$\Box$

\medskip
Up to now, from subsections {\bf 4.2} --- {\bf 4.4}, we have finally
established the Drinfel'd isomorphism in the two-parameter case.

\newpage
\vskip30pt \centerline{\bf ACKNOWLEDGMENT}

\vskip15pt Part of this work was done when Hu visited l'DMA, l'Ecole
Normale Sup\'erieure de Paris from October to November, 2004, the
Fachbereich Mathematik der Universit\"at Hamburg from November 2004
to February 2005,  as well as the ICTP (Trieste, Italy) from March
to August, 2006. He would like to express his deep thanks to ENS de
Paris, H. Strade and ICTP for the hospitalities and the supports
from ENS, DFG and ICTP. Authors are indebted to the referee for the
useful comments.

\bigskip

\bigskip

\bibliographystyle{amsalpha}

\end{document}